\documentclass[10pt,twoside,final]{article}
\usepackage{enumerate}

\usepackage{hyperref}

\usepackage{geometry}

\usepackage{amsmath, amssymb, amsthm}

\usepackage[dvipsnames]{xcolor}

\usepackage{graphicx}
\usepackage{float}
\usepackage{multirow}
\usepackage{subcaption}

\usepackage{wrapfig}

\newcommand{\abs}[1]{\lvert #1 \rvert}
\newcommand{\norm}[1]{\Vert #1 \Vert}

\renewcommand{\P}{\mathbb{P}}

\def\CFL{\text{CFL}}

\newcommand{\cT}{\mathcal{T}}

\newcommand{\xLeftrightarrow}[2][]{\ext@arrow 0099\Leftrightarrowfill@{#1}{#2}}
\newcommand{\xRightarrow}[2][]{\ext@arrow 0099\Rightarrowfill@{#1}{#2}}

\newcommand{\bsigma}{\boldsymbol{\sigma}}
\newcommand{\bxi}{\boldsymbol{\xi}}
\newcommand{\bxit}{\tilde{\boldsymbol{\xi}}}

\newcommand{\xx}{\mathbf{x}}
\newcommand{\vv}{\mathbf{v}}
\newcommand{\q}{\mathbf{q}}
\newcommand{\uu}{\mathbf{u}}
\newcommand{\f}{\mathbf{f}}
\newcommand{\g}{\mathbf{g}}
\newcommand{\h}{\mathbf{h}}
\newcommand{\Q}{\mathbf{Q}}
\newcommand{\K}{\mathbf{K}}

\newcommand{\N}{\mathcal{N}}
\newcommand{\Nst}{\mathcal{N}^{st}}
\newcommand{\bzero}{\mathbf{0}}

\usepackage{siunitx}
\let\leq\leqslant\let\geq\geqslant

\newcommand{\edmt}[3][\footnotesize]{
\par\noindent{#1\sffamily\bfseries #2~}{#1\rmfamily #3}\vskip.6\baselineskip}
\def\acknowl#1{\edmt{Acknowledgements}{#1}}
\def\funding#1{\edmt{Funding}{#1}}
\def\dataava#1{\edmt{Data Availability}{#1}}
\def\complia#1{\edmt[]{Compliance with Ethical Standards}{}\edmt{Conflict of Interest}{#1}}

\AtBeginDocument{%
\nocite{label}

}

\graphicspath{{figs/}}

\begin{document}

\title{Continuous 3D Finite Element Subgrid Basis Functions for Discontinuous Galerkin Methods on Polyhedral Meshes}

\author{
Sixtine Michel\thanks{De Vinci Higher Education, De Vinci Research Center, Paris, France.}
\and
Lorenzo Diazzi\thanks{CNR IMATI, Genova, Italy.}
\and
Walter Boscheri\thanks{Laboratoire de Mathématiques UMR 5127 CNRS, Université Savoie Mont Blanc, France} \thanks{
Department of Mathematics and Computer Science, University of Ferrara, Italy.}
}

\maketitle
\abstract{
We present a novel high-order accurate nodal discontinuous Galerkin (DG) method for solving nonlinear hyperbolic systems of partial differential equations (PDEs) on fully unstructured three-dimensional polyhedral meshes. A mesh generator is firstly discussed in detail, which ensures the generation of admissible control volumes. For the first time, we then extend the concept of agglomerated finite element (AFE) basis functions to polyhedral grids. In this context, the discrete solution is represented within each polyhedral element using piecewise continuous polynomials of degree $N$, defined on an internal tetrahedral subgrid. 
The AFE basis functions are therefore constructed by agglomerating standard finite element basis functions on each sub-tetrahedron of the computational cell. This allows for the precomputation of universal local matrices (mass and stiffness) on the reference element given by the unit tetrahedron, enabling a quadrature-free implementation that remains efficient even on highly irregular polyhedral meshes. 
High-order of accuracy in time is achieved using a local space-time Galerkin predictor as part of the ADER approach, applied independently within each polyhedral element. 
To ensure robustness in the presence of discontinuities such as shocks, an artificial viscosity limiter is embedded into the numerical scheme, allowing for controlled dissipation and stabilization without compromising the overall accuracy in smooth regions. 
To demonstrate the robustness and accuracy of the method, we validate it through different three-dimensional benchmark problems for the compressible Euler and Navier-Stokes equations. This work opens the door to the use of quadrature-free high-order one-step DG methods on general three-dimensional polyhedral meshes, providing increased flexibility in mesh generation and adaptation for complex geometries, while maintaining computational efficiency and high-order accuracy.
}

\noindent\textbf{Keywords:} Continuous finite element subgrid basis for DG schemes, high order quadrature-free ADER-DG schemes, unstructured polyhedral meshes, compressible Euler and Navier--Stokes equations.

\medskip

\noindent\textbf{MSC:} 65M60, 65M12, 76M10, 35L65

\section{Introduction} \label{sec:intro}

The study of physical phenomena such as fluid dynamics, wave propagation, and heat transfer is a central topic in applied mathematics and engineering. Over the past few decades, this area has benefited from advances in computer performance \cite{Zelkowitz2009AdvancesIC} and the development of sophisticated numerical methods \cite{bookMastorakis2009} to accurately capture the behavior of complex physical systems, thereby enabling numerical simulations. 
In particular, the time evolution of physical quantities such as mass, momentum, or energy in air or liquids can often be described by nonlinear systems of partial differential equations (PDEs), e.g. the Euler and Navier-Stokes equations govern the dynamics of compressible inviscid and viscous flows, respectively. When extended to multidimensional and realistic domains, the resolution of these systems becomes particularly challenging and requires efficient and accurate numerical strategies. Furthermore, the geometrical complexity of real-world domains often requires the use of unstructured and adaptive meshes (locally refined in regions of interest), and the nonlinear nature of the governing equations can lead to the formation of discontinuities in the solutions which pose challenges in the development of numerical schemes. 

Indeed, as hyperbolic PDEs are often characterized by solutions that develop discontinuities, even from smooth initial data, it is necessary to accurately capture these features without generating spurious oscillations or numerical instability. In continuous Galerkin (CG) methods, stabilizing the numerical scheme typically requires the addition of artificial diffusion terms \cite{KuzminStabilizationCG} (as also reported in \cite{Michel2023} and references therein). Discontinuous Galerkin (DG) methods use numerical fluxes and Riemann solvers at element interfaces, which are naturally suited for handling discontinuities \cite{ReedHill1973}. Nevertheless, in both cases, additional artificial viscosity techniques have to be introduced to further suppress non-physical oscillations near shocks, while maintaining high-order accuracy in smooth regions (e.g. shock capturing and subcell correction strategies  \cite{Gaburro2024}). 

Achieving high-order accuracy in time is also important, particularly for long-time unsteady simulations. Traditional time integration methods used to solve nonlinear hyperbolic conservation laws, such as Runge-Kutta (RK) or Strong Stability Preserving Runge-Kutta (SSPRK) schemes (see \cite{SHU1988439,Spiteri2006,Ruuth2006,Gottlieb2006OptimalST} and references therein), can be combined with high-order spatial discretizations, but may require an evaluation of the solution at each intermediate RK stage. One-step methods have then emerged from the ADER (Arbitrary high-order DERivatives) approach, originally proposed by Toro et al. \cite{Titarev2002,Titarev2005,TITAREV2005715,TORO2006150,BUSTO2016553}, which provides a fully-discrete one-step time integration technique that achieves arbitrary high-order accuracy in both space and time \cite{DUMBSER20088209}. 
This method relies on a local space-time Galerkin predictor, which evolves high-order polynomials in both space and time within each element. These polynomials are then used to evaluate numerical fluxes at cell interfaces. The ADER strategy is particularly well-suited for DG methods and enables high efficiency in large-scale parallel simulations. We refer also to \cite{inbookToro2024} and references therein for a broad overview of the ADER methodology combined with different spatial discretisations.

Although a lot of numerical schemes have been developed and extensively studied in the context of one- and two-dimensional problems, their extension to three-dimensional problems remains less developed. This is particularly due to the algorithmic complexity of handling unstructured three-dimensional meshes, the increased number of degrees of freedom to take into account, and the inherent complexity of high-order schemes which significantly increases the computational cost. As a result, there is a clear need of developing advanced numerical methods for fully three-dimensional simulations \cite{BOSCHERI2022127457}. 

A particularly elegant and efficient strategy to simultaneously handle geometric flexibility, high-order accuracy, and the presence of discontinuities has been recently introduced in \cite{BoscheriAFE2022} for two-dimensional problems. There, a discontinuous Galerkin (DG) method is combined for the first time with an agglomerated (continuous) finite element (AFE) framework, offering a unified solution to several challenges in the numerical resolution of PDE systems. In particular, classical finite element or finite volume methods are typically designed for meshes made of triangles or quadrilaterals in 2D, and tetrahedra or hexahedra in 3D. However, in real-world applications, such more regular meshes are often insufficient to capture complex geometries, and polytopal meshes become extremely appealing. 
Despite the advantages, the main problem related to polytopal tessellations lies in the fact that no reference element can be identified, thus preventing precomputation of mass and stiffness matrices in the preprocessing stage. Even the Virtual Element Method \cite{vem2,vem3,vem4,vem5} requires element-dependent projection matrices to be evaluated, despite its great flexibility in handling polytopal meshes, even in the non-conforming setting. The method introduced in \cite{BoscheriAFE2022} allows Voronoi meshes to embed a high-order DG method capable of operating on general unstructured grids, thus avoiding the limitations of classical element shapes. Then, in this direction, to address the use of polygonal or polyhedral meshes, the concept of Agglomerated Finite Elements (AFE) is introduced, which involves subdividing each polygonal cell into a set of triangular sub-elements, upon which standard finite element basis functions are defined. The global solution within each polygonal cell is then represented as piecewise continuous polynomial functions, defined over the subgrid. This process results in a piecewise continuous finite element basis inside each arbitrary polygonal element, while preserving the geometric flexibility of the original mesh, and keeping a discontinuous solution across inter-element boundaries.

In this work, we propose a novel high-order discontinuous Galerkin (DG) method on unstructured three-dimensional polyhedral meshes. We extend the concept of agglomerated finite element (AFE) basis, previously developed in \cite{BoscheriAFE2022} for two-dimensional polygonal grids, to the three-dimensional setting, enabling efficient and accurate simulations on general polyhedral elements. To stabilize the new DG schemes across shocks and discontinuities, we rely on an artificial viscosity limiter to handle shocks and discontinuities robustly.

The paper is organized as follows. In Section \ref{sec.PDE} we formulate the mathematical problem under consideration. Section \ref{sec:mesh} presents the three-dimensional mesh generation algorithm, highlighting the special ingredient of the AFE approach. Section \ref{sec:num_scheme} describes the numerical scheme, including the construction of the AFE basis functions and the ADER time integration method. Section \ref{sec:num_res} reports a series of three-dimensional benchmark tests to demonstrate the accuracy (through convergence tests) and robustness of the method. Finally, Section \ref{sec:conclu} provides some concluding remarks and proposes future perspectives.


\section{Governing equations} \label{sec.PDE}
We consider the compressible Navier-Stokes equations written in the following general formulation for nonlinear hyperbolic systems of conservation laws:
\begin{equation}
\frac{\partial \mathbf{Q}}{\partial t} + \nabla \cdot \mathbf{F}(\mathbf{Q}, \nabla \mathbf{Q}) = \mathbf{S}(\mathbf{Q}), \quad \mathbf{x} = (x,y,z) \in \Omega \subset \mathbb{R}^3, \quad t \in \mathbb{R}^+,
\label{eq:conservation_law}
\end{equation}
in the bounded domain $\Omega \subset \mathbb{R}^3$, where $\mathbf{Q} = \mathbf{Q}(\mathbf{x},t)$ is the vector of conserved variables, $\mathbf{F}(\mathbf{Q}, \nabla \mathbf{Q})$ represents the total flux tensor, including both inviscid and viscous contributions, and $\mathbf{S}(\mathbf{Q})$ is the source term, accounting for possible external forces. The vector of conserved variables is given by $\mathbf{Q} = (\rho,\rho \mathbf{v}, \rho E)^T$ with $\rho$ being the density, $\mathbf{v}=(u,v,w)$ the velocity vector and $\rho E$ the total energy. The nonlinear flux tensor reads
\begin{eqnarray}
\mathbf{F}(\mathbf{Q}, \nabla \mathbf{Q}) &=& (\mathbf{f}(\mathbf{Q}, \nabla \mathbf{Q}), \mathbf{g}(\mathbf{Q}, \nabla \mathbf{Q}), \mathbf{h}(\mathbf{Q}, \nabla \mathbf{Q})) \nonumber \\
&=&	\left( \begin{array}{c}  \rho \vv \\ \rho \left(\vv \otimes \vv \right) + \boldsymbol{\sigma}(\Q,\nabla \Q) \\ \mathbf{v} \cdot (\rho E \mathbf{I} + \boldsymbol{\sigma}(\Q,\nabla \mathbf{Q}) ) - \kappa \nabla T  \end{array} \right),
\end{eqnarray}
where $\mathbf{I}_{[3 \times 3]}$ is the identity matrix and $\bsigma$ is the viscous stress tensor given by the Stokes hypothesis \cite{Stokes1845InternalFriction} as 
\begin{equation}
\bsigma = \left( p + \frac{2}{3} \mu \, \nabla \cdot \mathbf{v} \right) \mathbf{I} - \mu \left( \nabla \mathbf{v} + \nabla \mathbf{v}^\top \right),
\label{eq:stress_tensor}
\end{equation}
with $\mu$ representing the dynamic viscosity coefficient. The temperature is denoted by $T$ and $\kappa$ is the heat conduction coefficient. The system is closed by the ideal gas equation of state:
\begin{equation}
p = (\gamma - 1) \left( \rho E - \frac{1}{2} \rho \|\vv \|^2 \right), \label{eq:EOS}
\end{equation}
where $\gamma$ is the ratio of specific heats. For the three-dimensional compressible Navier-Stokes equations, the convective eigenvalues $\boldsymbol{\lambda}$ and the viscous eigenvalues $\boldsymbol{\lambda}^v$ are given by \cite{DUMBSER201060}
\begin{equation}
\boldsymbol{\lambda} = (\abs{\vv} + c,\, \abs{\vv},\, \abs{\vv},\, \abs{\vv},\, \abs{\vv} - c), \qquad
\boldsymbol{\lambda}^v = \left( \frac{4 \mu}{3 \rho}, \frac{\gamma \mu}{\mathrm{Pr} \, \rho} \right),
\label{eq:eigenvalues3D}
\end{equation}
with $c^2 = \gamma R T$ being the speed of sound and $R$ the gas constant. The Prandtl number is given by $\mathrm{Pr}=\mu \gamma c_v / \kappa$.

\section{The three-dimensional mesh generation} \label{sec:mesh}

Polyhedral meshing algorithms use generic polyhedra to subdivide a given volume which defines the computational domain. Typically, tetrahedral and hexahedral meshing techniques are preferred as they lead to considerable simplifications in the algorithm design and implementation. Nonetheless, with respect to tetrahedral and hexahedral meshes, polyhedral meshes have the advantage of describing a computational domain using fewer elements but comparable accuracy. This fact is particularly significant when computing numerical solutions, as it leads to reduce computational costs.  Polyhedral meshing has received increasing attentions in the last decades due to the development of advanced numerical schemes in the finite volume and finite element context. Existing polyhedral meshing techniques can be grouped depending on the way the domain subdivision is approached. Methods employing Voronoi tessellation \cite{Gar2014,LEE2015,YAN2013} compute the dual of a Delaunay tetrahedral mesh to obtain a polyhedral mesh formed by Voronoi cells. Agglomeration strategies \cite{ANTONIETTI2025,CONTRERAS2014} essentially allow to define polyhedral cells by merging the elements of an initial tetrahedral mesh, possibly guided by some quality criteria. Another approach consists of using a binary space partition to split the convex hull defined by an input triangle soup into a constrained mesh of convex polyhedra \cite{DIAZZI2021}.

In this work our objective is to create a \textit{constrained} polyhedral mesh $\big\{P_i\big\}_{i=1}^{N_e}$ of the computational domain $\Omega$, where the boundary of $\Omega$ is represented as a union of facets of mesh polyhedra $P_i$ with a total number $N_e$ of cells. Furthermore, it is required that each polyhedron admits a tetrahedrization and the collection of all the tetrahedra $T_i$ forms a valid sub-mesh $\big\{T_i\big\}_{i=1}^{N_t}$ composed of $N_t$ cells. To discretize each three-dimensional computational domain, we resort to a polyhedral meshing algorithm based on local operations that split and agglomerate elements of a starting constrained tetrahedral mesh $\big\{ \Gamma_i \big\}_{i=1}^{N_\gamma}$. 

\subsection{Local polyhedral replacement (LPR)} \label{sec:loc_poly_replacement}
The key idea of our approach is to combine two local operations to replace the set of tetrahedra incident at an internal vertex $v$ with a collection of polyhedra. A mesh vertex $v$ will be called a \textit{replacement-center} if the following conditions hold: (i) $v$ is neither on the domain boundary $\partial\Omega$ nor on the mesh convex hull, (ii) all the mesh elements incident at $v$ are tetrahedra.

\begin{figure}[!htbp]
    \centering
    \includegraphics[width=\textwidth]{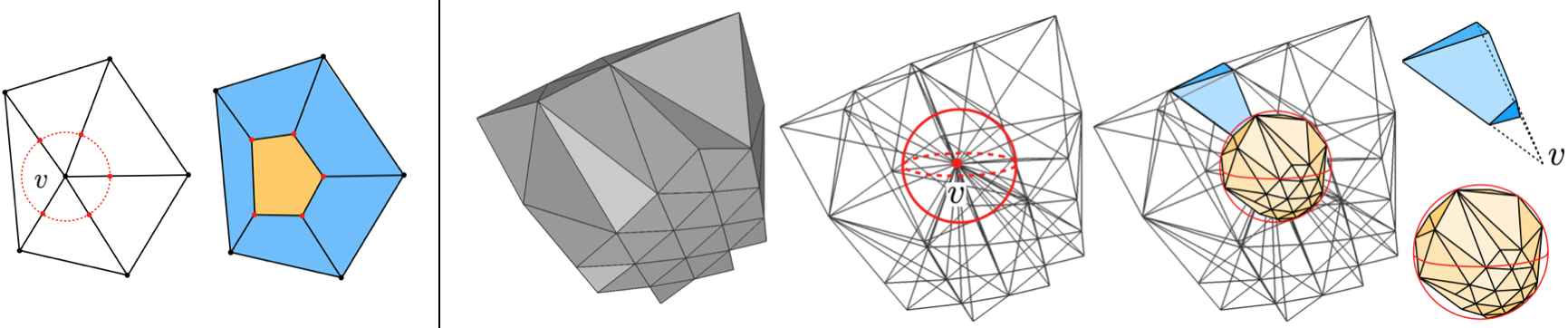}
    \caption{Two examples of LPR. Left: two-dimensional case showing how the local polyhedral replacement works. Right: three-dimensional case where we show (from left to right) the front view of the tetrahedra incident at a replacement-center $v$, the edges of these tetrahedra along with the sphere (red) picturing the distance $d$ from $v$ at which the split points are placed, and finally one of the truncated tetrahedra (blue) and the central polyhedron (yellow).}
    \label{fig:lpr_2D_3D_ex}
\end{figure}

\noindent
\begin{minipage}{0.8\linewidth}
\raggedright
Given $v$, we denote by $\Gamma_v$ the set of its $n$ incident tetrahedra, and by $E_v$ the set of all the tetrahedra edges incident at $v$. The first operation consists of splitting each edge in $E_v$ by inserting on it a new vertex at distance $d$ from $v$, where $d= \frac{1}{2} \displaystyle{ \min_{e\in E_v}} \norm{e}$ and $\|.\|$ is the Euclidean norm. Each tetrahedron in $\Gamma_v$ is then split by inserting the triangle defined by the three split points placed on its edges incident at $v$ (see the schematic on the right) producing two polyhedral cells: a new smaller tetrahedron $T_i$ incident at $v$ and a truncated tetrahedron $S_i$. 
\end{minipage}
\hfill
\begin{minipage}{0.2\linewidth}
\includegraphics[width=\linewidth]{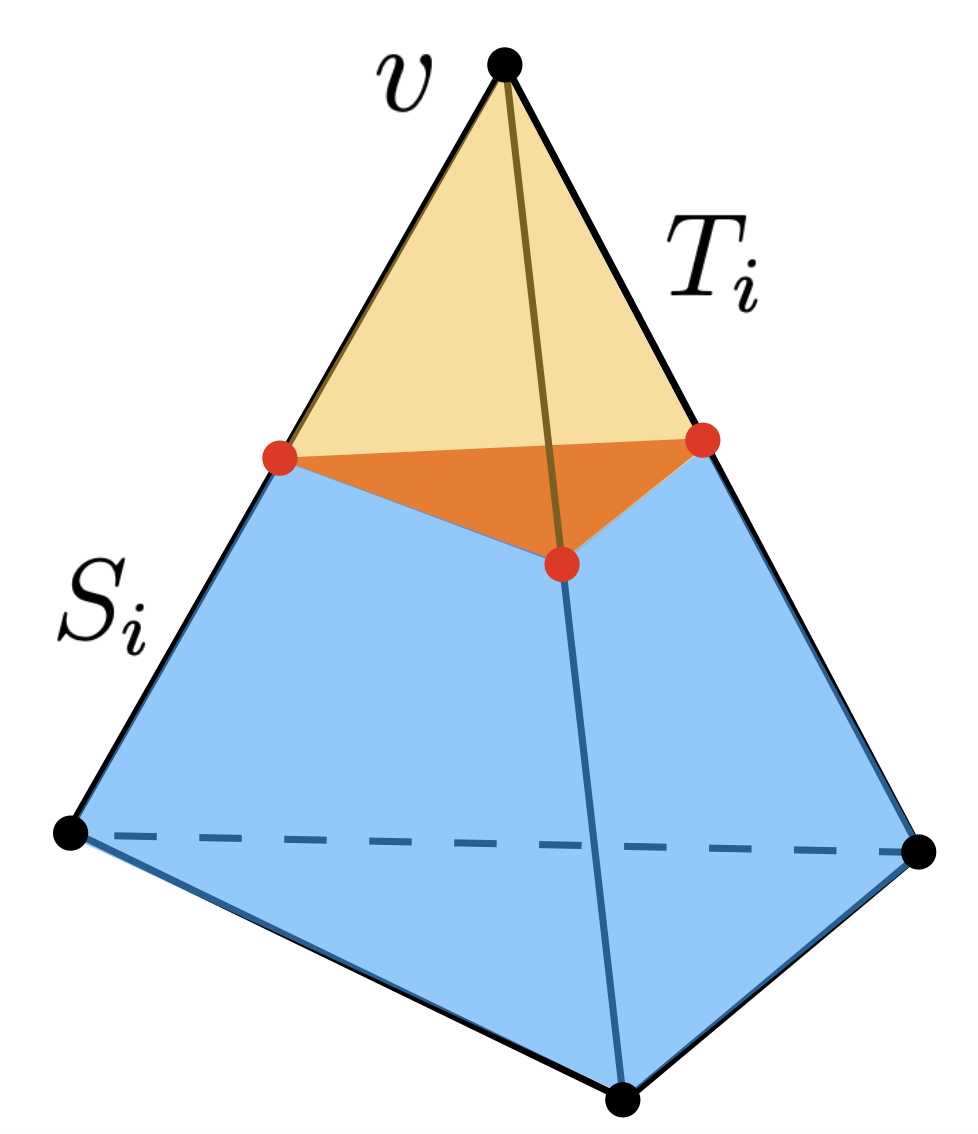}
\end{minipage}
Each of the truncated tetrahedra is a pentahedron with two triangular faces and three quadrilateral faces. The second operation relies on melting all the new tetrahedra $T_i$ in a unique \textit{central-polyhedron} $P$. All the triangles incident at $v$ are removed and the $n$ remaining triangles opposite to $v$ become the faces of $P$. A two-dimensional example, to give a simple and clear idea of the procedure, and also a well-suited three-dimensional example are shown in Fig. \ref{fig:lpr_2D_3D_ex}.

\begin{figure}[!htbp]
    \centering
    \includegraphics[width=0.5\textwidth]{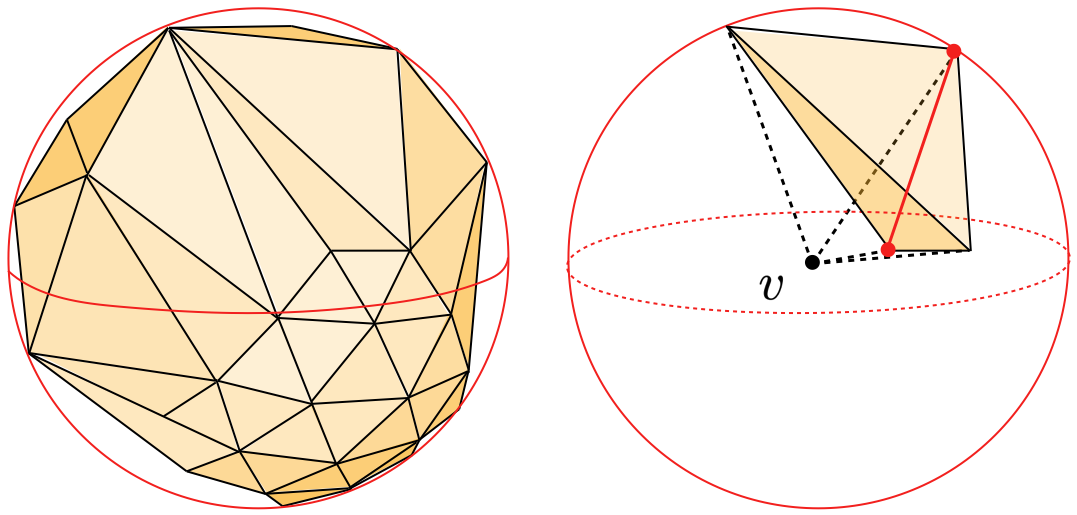}
    \caption{This example shows a concave central polyhedron. On the left only the front polyhedral triangular facets are visible. All the polyhedron triangular faces have their vertices on the red sphere whose center is the vertex $v$. The polyhedron is concave as considering the two adjacent boundary triangles on the right, the red line connecting the vertices opposite to the shared (black) edge does not intersect either the interior or the boundary of the polyhedron.}
    \label{fig:non_convex_ex}
\end{figure}

\begin{figure}[!htbp]
    \centering
    \includegraphics[width=\linewidth]{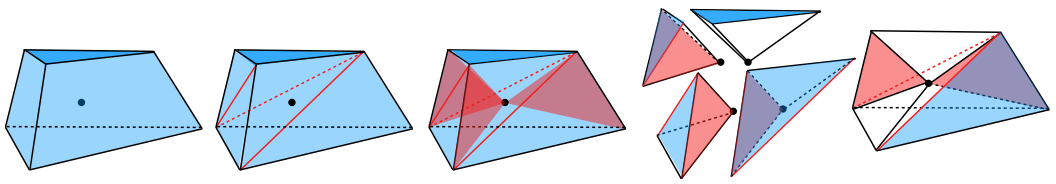}
    \caption{Tetrahedrization of a truncated tetrahedron $S_i$. On the left $S_i$ is shown (blue) with its barycenter. The quadrilateral faces are triangulated by fixing three diagonals (red segments), while the interior is divided by adding $12$ triangles connecting each edge or diagonal to the barycenter (the triangles incident at edges connecting the two triangular faces of $S_i$ are highlighted in red). On the right, four tetrahedra are detached to further clarify the three-dimensional structure.}
    \label{fig:trucated_tet_tetrahedrization}
\end{figure}

LPR at $v$ modifies the mesh only with respect to $\Gamma_v$, keeping all the other mesh elements untouched. In particular, the boundary of $\Gamma_v$ is preserved, thus ensuring that the mesh connectivity remains valid after the element substitution. Furthermore, this procedure increments the total number of mesh elements of exactly one unit, the $n$ tetrahedra in $\Gamma_v$ are replaced by $n$ truncated tetrahedra plus one \textit{central-polyhedron}. Since all the polyhedra created by LPR are tetrahedrizable, it is ensured the existence of a tetrahedral sub-mesh, which is crucial for the definition of the agglomerated finite element basis functions used by the discontinuous Galerkin scheme. The \textit{central-polyhedron} $P$, which may not be convex (see Fig. \ref{fig:non_convex_ex}), can be subdivided into tetrahedra by connecting each of its $n$ triangular faces with its circumcenter $v$. 
Any of the truncated tetrahedra $S_i$ is a convex polyhedron and therefore tetrahedrizable. Among all the possible ways in which $S_i$ may be subdivided, we choose to triangulate each of its $3$ quadrilateral faces by appropriately fixing a diagonal and subsequently using the barycenter of $S_i$ as apex to create a tetrahedron by connecting it with each of the $8$ triangular faces bounding $S_i$ (see Fig. \ref{fig:trucated_tet_tetrahedrization}). Since a quadrilateral face $f$ shared between two truncated tetrahedra $S_i$ and $S_i'$ has two diagonals, it is important to tetrahedralize $S_i$ and $S_i'$ in such a way that the same diagonal is fixed on $f$. Doing so each truncated tetrahedron becomes an octahedron with only triangular faces.

Summarizing, LPR replaces the $n$ tetrahedra in $\Gamma_v$ with a set of polyhedra $P_v$ made of $n$ octahedra plus one $n$-faces polyhedron, such that $P_v$ admits a tetrahedrization conforming to its polyhedral cells counting $9n$ new tetrahedra.

\subsection{Polyhedral mesh and tetrahedral subgrid} \label{sec:poly_mesh_tet_subgrid}
Given the computational domain $\Omega \subset \mathbb{R}^3$, the generation of the polyhedral mesh consists of two main phases. First, a \textit{constrained Delaunay tetrahedrization} (CDT)  $\big\{ \Gamma_i \big\}_{i=1}^{N_\gamma}$ with respect to the boundary of $\Omega$, i.e. $\partial \Omega$, is built using the efficient algorithm proposed in \cite{DIAZZI2023}. Successively, LPR is performed at each replacement-center to create a polyhedral mesh $\big\{ P_i \big\}_{i=1}^{N_e}$.

\begin{figure}[!hbtp]
    \centering
    \includegraphics[width=\linewidth]{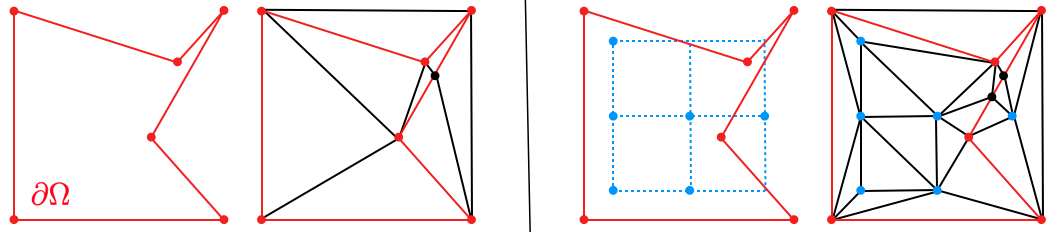}
    \caption{This two-dimensional example depicts a comprehensive picture of how grid vertices (blue) can be added to the computational domain to introduce valid replacement-centers. Being in two space dimensions, the domain is bounded by segments (red) rather than triangles. On the left, the CDT of the domain is shown in which a Steiner point (black vertex) has been inserted. Note that all the vertices are input vertices (red) or Steiner points (black) belonging to $\partial\Omega$. On the right there is the CDT of the enriched domain shown. Only grid based vertices (blue) that do not overlap $\partial\Omega$ have been added.}
    \label{fig:input_enrich_2D}
\end{figure}

The input of the CDT algorithm is a manifold, but not necessarily connected, triangulated surface embedded in the three-dimensional space, the output is a tetrahedrization of the input's convex hull such that each input triangle is represented as the union of some output tetrahedra facets. Typically, to build the CDT a certain number of input segments are split at appropriate points which are denoted as \textit{Steiner points}. All the CDT mesh vertices are input points or Steiner points, both belonging to $\partial\Omega$ and then not suitable to be used as replacement-centers. To fix this issue we choose a discretization step size $h$ (much smaller than the size of the domain) to generate a virtual uniform cubic grid fitting $\Omega$ without overlapping its boundary. These grid vertices and the boundary $\partial\Omega$ are joined to obtain a valid input for the CDT algorithm (see Fig. \ref{fig:input_enrich_2D}) such that LPR can be performed on the resulting mesh $\big\{\Gamma_i\big\}_{i=1}^{N_\gamma}$.

At the beginning, all replacement centers in $\big\{\Gamma_i\big\}_{i=1}^{N_\gamma}$ are collected in a list whose first vertex at each iteration is used to perform local polyhedral replacement. Once the list becomes empty, the algorithm terminates and the initial tetrahedral mesh $\big\{\Gamma_i\big\}_{i=1}^{N_\gamma}$ has been turned in to the final polyhedral mesh $\big\{P_i\big\}_{i=1}^{N_e}$. Note that, if $m$ local polyhedral replacement steps have been performed, the final mesh $\big\{P_i\big\}_{i=1}^{N_e}$ counts $N_e = N_\gamma + m$ polyhedral cells, while the tetrahedral sub-grid $\big\{T_i\big\}_{i=1}^{N_t}$ counts $N_t = N_\gamma + 8 ( n_1 + \dots + n_m ) $ elements where $n_i$ is the number of tetrahedra incident at replacement-center $v_i$ before substitution.

\begin{figure}[!hbtp]
	\centering
	\includegraphics[width=0.9\linewidth]{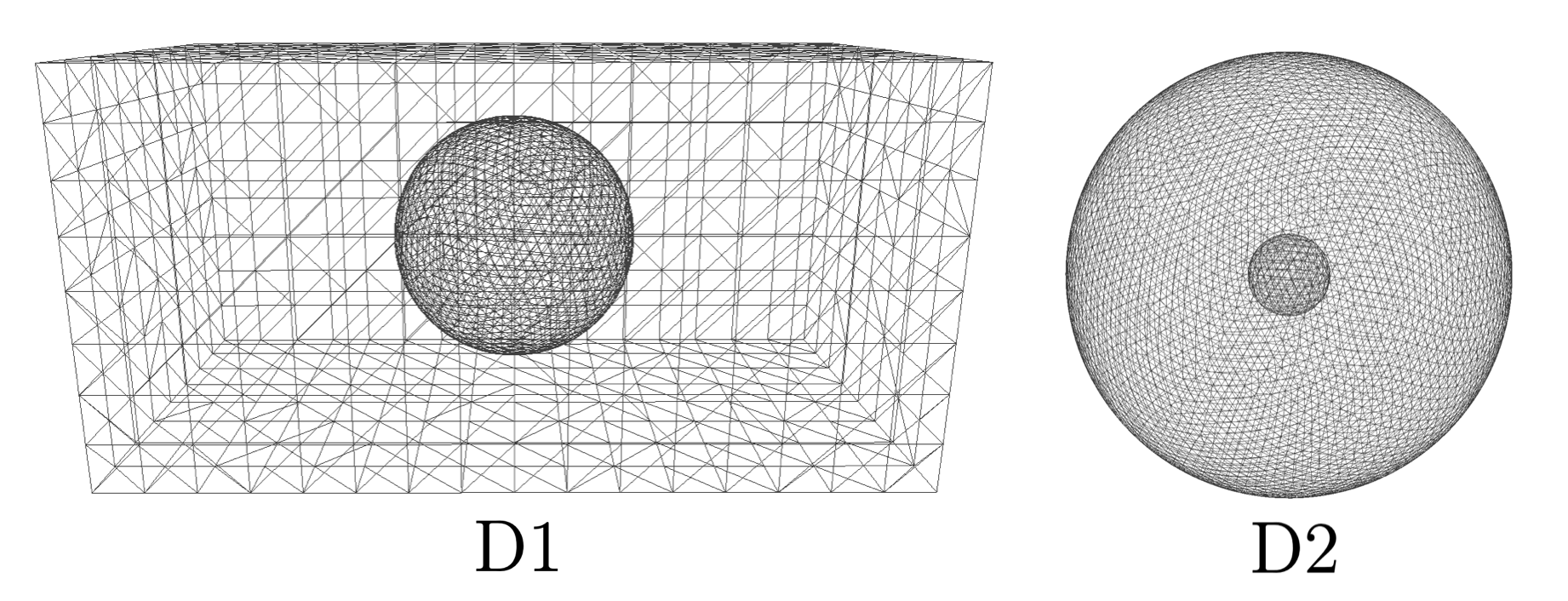}
	\caption{The picture shows the two domains used to collect the data reported in Table \ref{tab:mesh_gen_compare}. The domain referred as D1, on the left, is made by two concentric surfaces. The outer one is the triangulated boundary of a parallelepiped centered at the origin and whose sides are aligned with the Cartesian axes. The length of the sides along the $x-$direction is $6$, while the length of the others is $3$. The triangulated surface inside the parallelepiped has all its vertices on the unit sphere. On the right the domain referred as D2 is made by two triangulated surfaces whose common center is the origin and having their vertices on the spheres of radius $1$ and $0.2$ respectively. }
	\label{fig:two_domains}
\end{figure}

To provide an overview of the structure of the polyhedral meshes generated by LPR, we consider the computational domain depicted in Fig. \ref{fig:two_domains}. Some significant data collected from these examples are reported in Table \ref{tab:mesh_gen_compare}. We compare LPR with other two, more trivial, polyhedral mesh generation techniques based on a straightforward vertex or edge agglomeration (briefly described in the caption of Table \ref{tab:mesh_gen_compare}). As one could expect, in these examples LPR generates the lowest percentage of tetrahedra, hence maximizing the number of polyhedra in the computational grid.
\begin{table}[!htbp]
    \centering
    \renewcommand{\arraystretch}{1.5}
    \begin{tabular}{|c| c | c | c | c | c | c | c | c | c |}
        \hline
        \multicolumn{2}{|c|}{Dom.} & $\frac{\sqrt[3]{Vol}}{h}$ & $|V|$ & $n_e$ & $n_4$ (\%) & $n_6$ (\%) & $n_8$ (\%) & $n_{10}$ (\%) & $n_{>10}$ (\%) \\ \hline \hline
        \multirow{6}{*}{\rotatebox{90}{LPR}}
        & D1 & 10 & 4680 & 12253 & 67.60 & 0.00 & 31.53 & 0.008 & 0.85 \\ \cline{2-10}
        & D1 & 20 & 24624 & 54504 & 49.55 & 0.00 & 48.41 & 0.007 & 2.03\\ \cline{2-10}
        & D1 & 40 & 191552 & 397922 & 41.25 & 0.00 & 56.28 & 0.005 & 2.46\\ \cline{2-10}
        & D2 & 10 & 3805 & 7767 & 56.73 & 0.00 & 41.84 & 0.05 & 1.38\\ \cline{2-10}
        & D2 & 20 & 21888 & 45443 & 41.57 & 0.00 & 56.01 & 0.007 & 2.41\\ \cline{2-10}
        & D2 & 40 & 187869 & 390479 & 39.18 & 0.00 & 58.22 & 0.01 & 2.59\\ \hline\hline
        \multirow{6}{*}{\rotatebox{90}{VT}}
        & D1 & 10 & 2585 & 7871 & 97.80 & 0.04 & 0.04 & 0.14 & 1.98\\ \cline{2-10}
        & D1 & 20 & 9628 & 27128 & 95.69 & 0.01 & 0.007 & 0.08 & 4.20\\ \cline{2-10}
        & D1 & 40 & 62990 & 184094 & 95.29 & 0.002 & 0.005 & 0.03 & 4.67\\ \cline{2-10}
        & D2 & 10 & 1995 & 3290 & 86.84 & 0.00 & 0.00 & 0.30 & 12.86\\ \cline{2-10}
        & D2 & 20 & 7167 & 19116 & 93.66 & 0.00 & 0.00 & 0.00 & 6.34\\ \cline{2-10}
        & D2 & 40 & 55876 & 172991 & 95.16 & 0.00 & 0.01 & 0.01 & 4.82\\ \hline\hline
        \multirow{6}{*}{\rotatebox{90}{ET}}
        & D1 & 10 & 2643 & 4715 & 59.15 & 5.00 & 9.86 & 8.02 & 17.96\\ \cline{2-10}
        & D1 & 20 & 10317 & 20360 & 57.47 & 3.02 & 12.23 & 12.47 & 14.81\\ \cline{2-10}
        & D1 & 40 & 69754 & 141839 & 56.66 & 1.55 & 11.54 & 13.39 & 16.85\\ \cline{2-10}
        & D2 & 10 & 2069 & 3105 & 42.93 & 25.18 & 14.40 & 7.47 & 10.02\\ \cline{2-10}
        & D2 & 20 & 8063 & 16181 & 54.84 & 4.19 & 11.95 & 12.46 & 16.56\\ \cline{2-10}
        & D2 & 40 & 64050 & 137622 & 57.67 & 1.065 & 9.95 & 13.93 & 17.38\\ \hline
    \end{tabular}
    \caption{Data from polyhedral meshes generated with different types of tetrahedra agglomeration strategies on different domains. VT and ET denote two agglomeration methods more trivial than LPR. Similarly to LPR, VT creates polyhedra by agglomerating tetrahedra incident at a vertex, but without splitting incident edges. ET creates polyhedra by merging all tetrahedra incident at an edge. Both VT and ET agglomeration strategies can be performed only at vertices/edges having only tetrahedra as incident elements. Domain D1 and D2 are described in Fig. \ref{fig:two_domains}. Each domain is enriched with vertices taken from a uniform grid of step-size $h$, centered at the origin and filling the domain. The second column reports the ratio between the cubic root of the domain volume and $h$. $|V|$ is the number of vertices of the output mesh. $n_e$ is the number of polyhedra of the output mesh. Columns tagged as $n_i$ (\%) report the percentage of polyhedra with $i$ facets in the output mesh ($i=4$ means tetrahedra, $i=6$ hexahedra, and so on).}
    \label{tab:mesh_gen_compare}
\end{table}

\subsection{Robust and efficient implementation}
The polyhedral meshing algorithm described in Section \ref{sec:poly_mesh_tet_subgrid} has been developed on the existing implementation of the CDT algorithm \cite{DIAZZI2023} (referred as $cdt\_2023$ hereafter). 
$cdt\_2023$ combines robustness and efficiency thanks to the usage of the so called \textit{implicit points} and \textit{indirect predicates} \cite{ATTENE2020}. The first part of our algorithm consists of computing the CDT of the computational domain $\partial\Omega$ which is natively supported by $cdt\_2023$, thus only functionalities dealing with LPR have been added, taking care of preserving $cdt\_2023$ robustness and efficiency. In order to build the CDT, the insertion of one or more Steiner points on input edges is required. While input points are represented using standard floating point coordinates, in $cdt\_2023$ Steiner points are encoded through a particular kind of implicit points named \textit{linear combination implicit points} (LNCs). As suggested by the name, a LNC defines a new point as the unevaluated linear combination $q_1 + k q_2$ with $k \in [0,1]$ of two input points $q_1$ and $q_2$. 
Employing a LNC ensures that the Steiner point is placed exactly on the edge it splits. On the contrary, rounding errors which may occur using floating point coordinates representation, may slightly move the Steiner point out from the edge, introducing invalid mesh configurations. 
To achieve robustness along the whole pipeline, LNCs are exploited during LPR to represent split points. \textit{Implicit point} representation is not supported by the numerical solver, thus, once the polyhedral mesh has been created, all LNC points are approximated to floating point coordinates. A post-processing mesh optimization feature is available in $cdt\_2023$ to improve mesh quality which, in practice, prevents the creation of invalid configurations during final coordinate rounding.

The elapsed time required by the algorithm $cdt\_2023$ strongly depends on the number of input vertices and triangles, and also on the way they are arranged i.e. on the geometrical structure of the input surface. Efficiency has been empirically proved in \cite{DIAZZI2023} by both showing that the measured elapsed time has $O(n\log n)$ behavior with respect to the number $n$ of required Steiner points and that performances are comparable with the fastest (but not robust) state of the art tools. Concerning the generation of the polyhedral mesh, we observe that the computational cost depends on the number $m$ of LPR steps performed. Each step involves the tetrahedra incident at each replacement-center $v$. Considering that all the vertices edge-connected to $v$ cannot be replacement-centers in turn, we argue that the computational time required by all the $m$ steps of LPR is linearly related to the number of tetrahedra in the CDT, thus ensuring the efficiency of the polyhedral generation procedure.

\section{Numerical scheme} \label{sec:num_scheme}

We consider the computational domain $\Omega \subset \mathbb{R}^3$ discretized into a set of non-overlapping polyhedral elements $\{P_i\}_{i=1}^{N_e}$. The goal is to construct a high-order accurate fully-discrete one-step discontinuous Galerkin (DG) method on these general polyhedral meshes. In the sequel, we first detail the new ansatz for the agglomerated finite element basis functions, then we describe the local space-time Galerkin predictor, and finally we present the fully-discrete corrector based on the DG formulation. 

\subsection{Agglomerated finite element nodal basis functions}
\begin{figure}[htbp]
    \centering
    \includegraphics[trim=2cm 1.75cm 1.75cm 2cm, clip, width=0.4\textwidth]{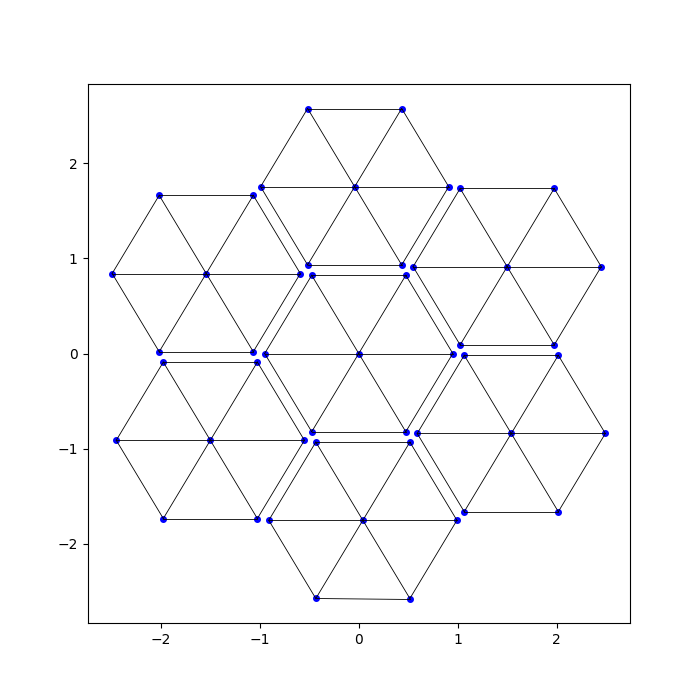} 
    \includegraphics[trim=2cm 1.75cm 1.75cm 2cm, clip, width=0.4\textwidth]{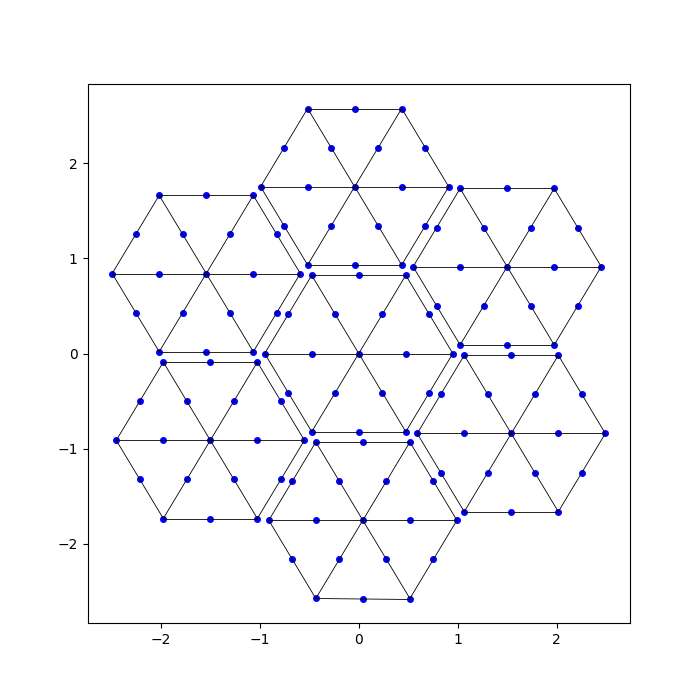} 
    \caption{{Representation of a two-dimensional AFE-DG discretization, with second (left) and third (right) order spatial accuracy. The DG space is discontinuous across polygonal elements and the CG space is piecewise continuous within each polygon on the triangular subgrid.}}
    \label{fig:AFE_2D}
\end{figure}

Unlike classical DG schemes that employ modal polynomial bases directly defined on the polyhedral cells, our approach uses nodal polynomial basis functions, based on the \textit{Agglomerated Finite Element} (AFE) basis ansatz proposed in \cite{BoscheriAFE2022}, which leverages a continuous Galerkin (CG) finite element space defined on a tetrahedral subgrid within each polyhedral element. A representation of the AFE-DG discretization is proposed in Fig. \ref{fig:AFE_2D} for a two-dimensional second and third order spatial discretization. In particular, we can see that the DG space is discontinuous across polygonal elements while the CG space remains piecewise continuous within each control volume by means of a triangular subgrid. 
This method is here extended to high-order nodal polynomial approximations on arbitrarily shaped polyhedra in three space dimensions, while preserving desirable properties such as locality, efficiency, and flexibility. 
The AFE framework allows the use of reference element formulations enabling the pre-computation of local matrices (mass, stiffness, flux) on a universal reference tetrahedron $\hat{T}$. This contributes to the development of a quadrature-free implementation, which remarkably reduces the computational cost without sacrificing accuracy.

\begin{figure}[!htbp]
	\centering
	\includegraphics[trim=6.5cm 8cm 6cm 8cm, clip, width=0.24\textwidth]{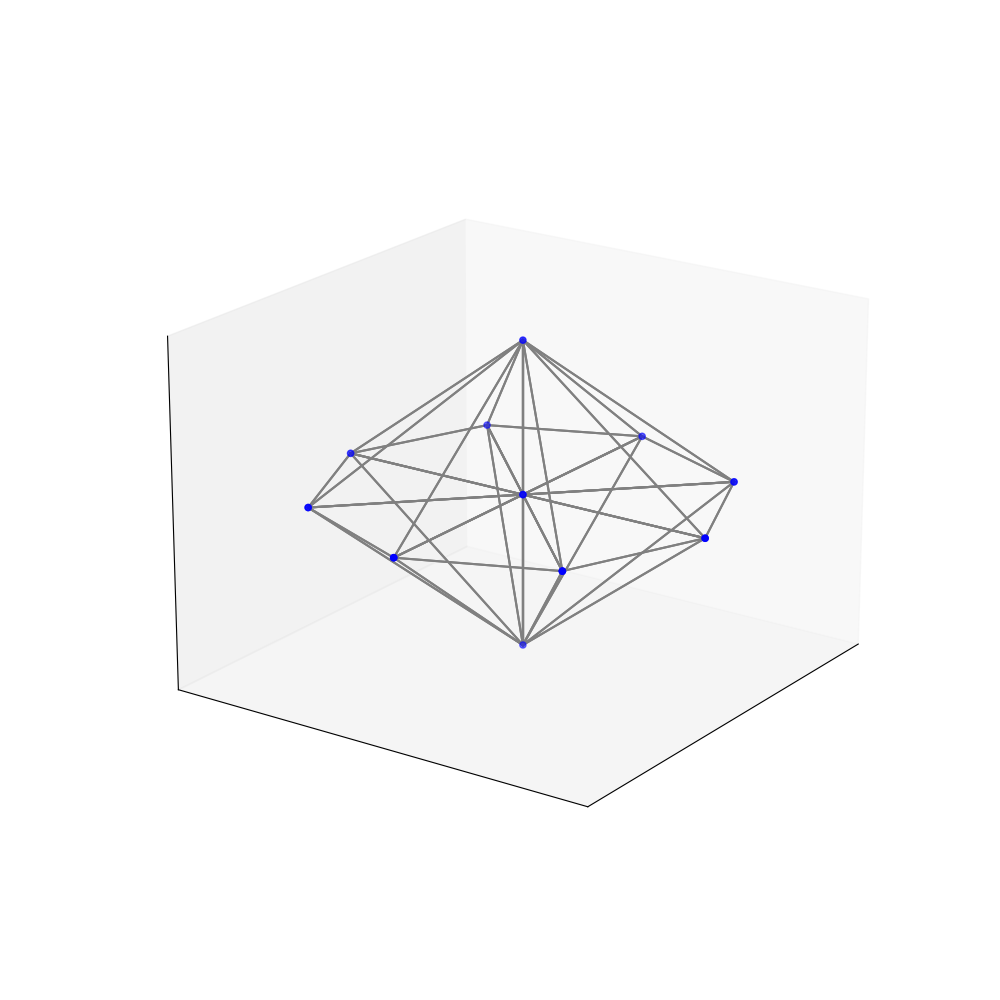}
	\includegraphics[trim=6.5cm 8cm 6cm 8cm, clip, width=0.24\textwidth]{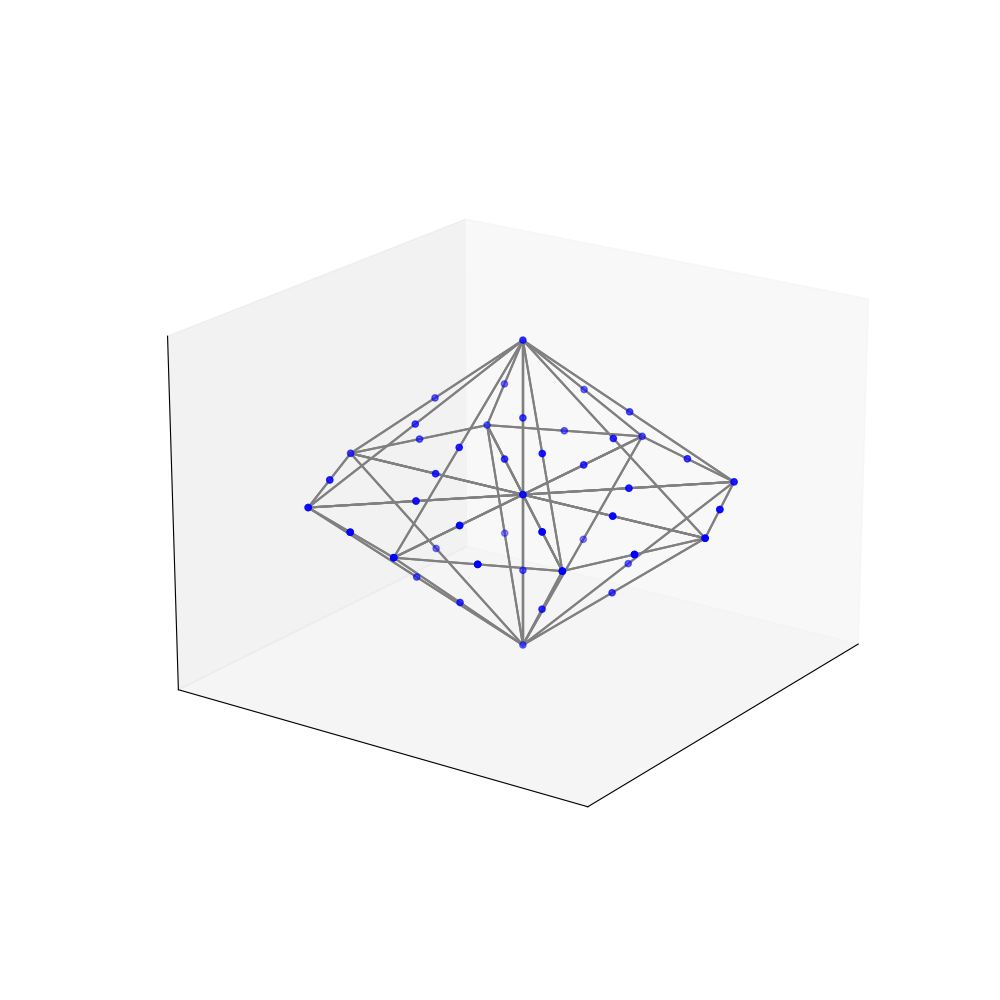}
	\includegraphics[trim=6.5cm 8cm 6cm 8cm, clip, width=0.24\textwidth]{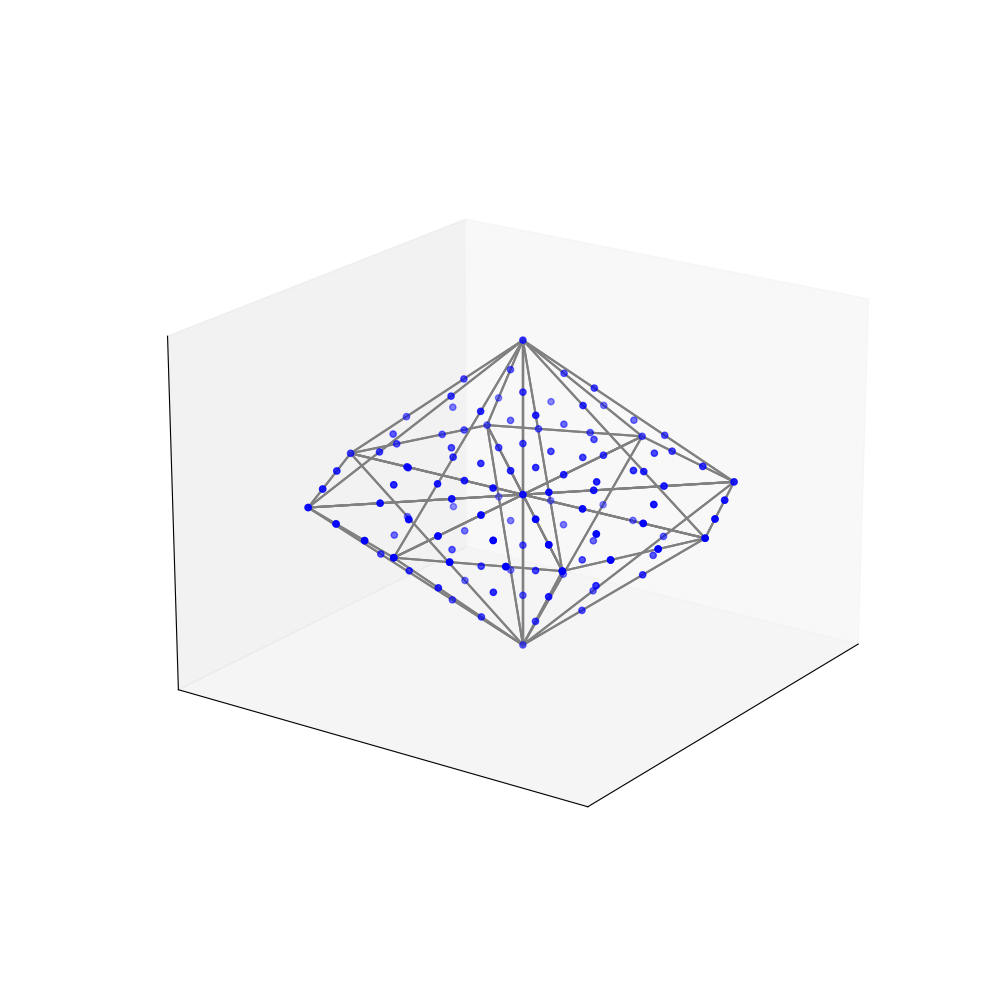}
	\includegraphics[trim=6.5cm 8cm 6cm 8cm, clip, width=0.24\textwidth]{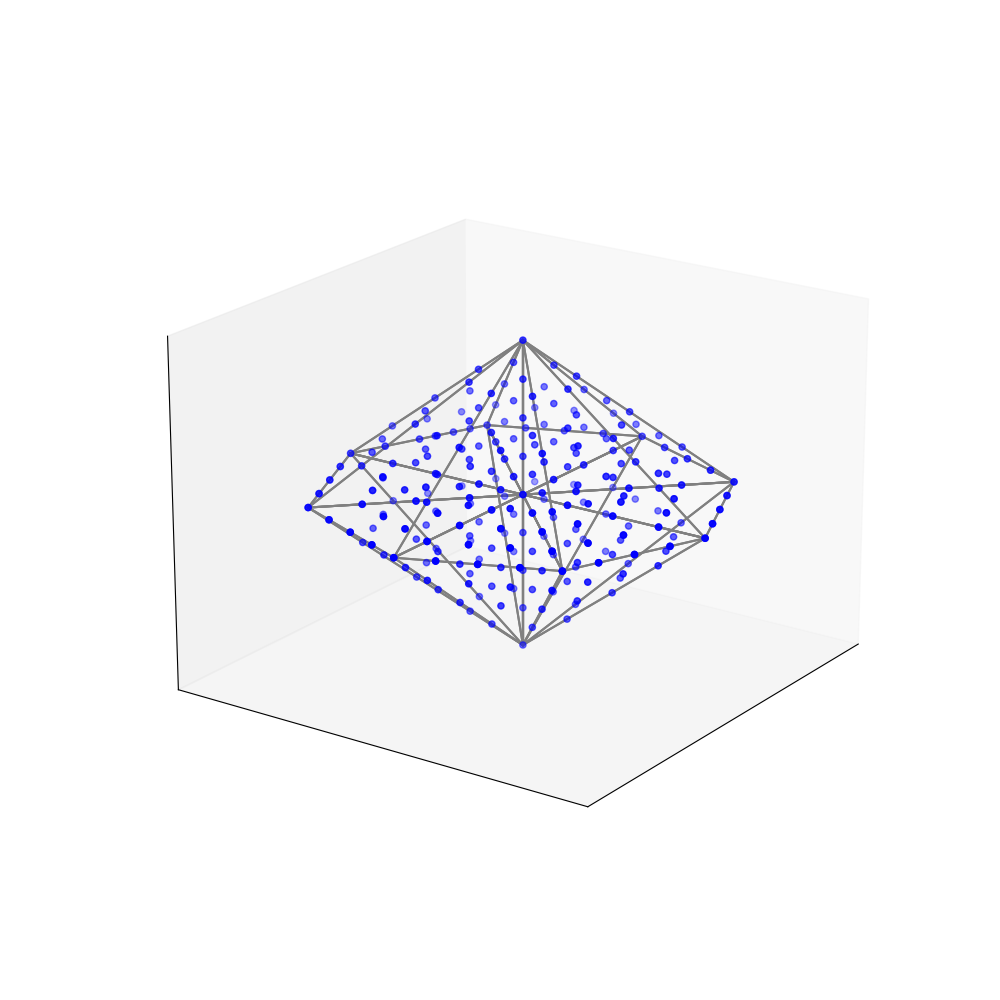}
    \caption{Example of a polyhedral element, divided into a set of tetrahedral elements with second- to fifth-order spatial discretizations, from left to right (i.e. the spatial discretization degree $N=\{1,2,3,4\}$).}
	\label{fig:AFE_all}
\end{figure}

Each polyhedral element $P_i$ is subdivided into a total number of $N_i$ tetrahedra $\{T_{ik}\}_{k=1}^{N_{i}}$, yielding a conforming subgrid within $P_i$. The subdivision can be performed by standard mesh generation algorithms, ensuring that the union of tetrahedra exactly covers $P_i$, and that tetrahedra meet only at faces, edges, or vertices, without any overlapping. An example of a tetrahedral subdivision is given in Fig. \ref{fig:AFE_all} for different spatial discretization orders.
On each tetrahedral sub-element $T_{ik}$, we denote by $N$ its spatial discretization degree, which leads to the definition of a $\mathbb{P}_N$ element. We consider the classical Lagrange finite element basis functions of degree $N$, denoted by $\{\phi^{(i,k)}_j(\mathbf{x})\}_{j=1}^{\mathcal{N}_p}$, where $\mathcal{N}_p = \frac{(N+1)(N+2)(N+3)}{6}$ is the number of basis functions for tetrahedral elements of degree $N$. These nodal basis functions satisfy the interpolation property, that is $\phi^{(i,k)}_j(\mathbf{x}^{(i,k)}_m) = \delta_{jm}$, where $\mathbf{x}^{(i,k)}_m$ are the physical coordinates of the interpolation nodes mapped from the reference element $\hat{T}$ to the physical tetrahedron $T_{ik}$. This property ensures that the degrees of freedom correspond directly to function values at these nodes, facilitating interpolation and projection operations. By employing a reference element mapping, all computationally intensive operations such as mass and stiffness matrix assembly or evaluation of gradients are performed once and for all on the reference tetrahedron $\hat{T}$. The transformation to the physical element is then handled via standard affine mappings, allowing for efficient and uniform implementation. 
The definition of these basis functions on the reference tetrahedral element for $N=\{1,2,3\}$ are fully detailed in Appendix \ref{app:basis_functions}. The discrete solution inside the whole polyhedral element $P_i$ is then represented as a piecewise continuous polynomial function, constructed by agglomerating the finite element basis functions defined on each tetrahedron of the internal subgrid associated to $P_i$:
\begin{equation}
\uu_h(\mathbf{x}, t) = \sum_{\ell=1}^{\N_{i}} \phi_\ell(\xx) \, \hat{\uu}_{\ell,i}(t) = \phi_\ell(\xx) \, \hat{\uu}_{\ell,i}(t), \quad \mathbf{x} \in P_i, \label{eq:interp_u}
\end{equation}
where $\hat{\uu}_{\ell,i}(t)$ are the time-dependent degrees of freedom (DOFs) which are given by the union, \textit{without repetition}, of all the degrees of freedom $\N_{i}$ contained in the subgrid made of $N_i$ tetrahedra for the cell $P_i$. Due to the continuity of the finite element basis within each $P_i$, the representation $\uu_h(\mathbf{x}, t)$ is piecewise continuous inside the element but discontinuous across inter-element boundaries. This discontinuity naturally fits the DG framework, where inter-element coupling is handled by numerical fluxes. 

\subsection{Time stepping discretization}
The time domain $[0,t_f]$ is partitioned into discrete intervals such that $t \in [t^n;t^{n+1}]$. A local parametrization of time within each interval is introduced as
\begin{equation}
	t = t^n + \tau \Delta t, \qquad \tau \in [0;1],
	\label{eqn.time_map}
\end{equation}
where $t^n$ denotes the current time level and $\Delta t := t^{n+1}-t^n$ is the time step size. The variable $\tau$ represents a reference time coordinate that maps the physical time interval onto the unit interval $[0,1]$ according to~\eqref{eqn.time_map}. 

The time step $\Delta t$ is chosen according to a CFL-type stability constraint for explicit DG schemes, given by
\begin{equation}
	\Delta t \leq \CFL \, \frac{\min \limits_{\Omega} h}{(2N+1) \max \limits_{\Omega} \left( \abs{\boldsymbol{\lambda}} + 2 \abs{\boldsymbol{\lambda}^v} \frac{2N+1}{h}\right)},
	\label{eqn.timestep}
\end{equation}
where $h$ is a characteristic length of element $P_i$ (e.g., the diameter of the polyhedral element) and $N$ is the polynomial degree of the approximation. The factor $(2N+1)$ arises from stability considerations specific to high-order DG methods \cite{CockburnKarniadakisShu2000,CockburnShu2001}. 

\subsection{Local space-time Galerkin predictor}
To achieve high-order accuracy in time, we employ the ADER (Arbitrary high-order DERivatives) approach firstly introduced in \cite{Titarev2002,Titarev2005}, and then generalized to deal with both finite volume and discontinuous Galerkin discretizations in \cite{DUMBSER20088209}. The ADER method is based on a local space-time Galerkin predictor $\mathbf{q}_h(\mathbf{x}, t)$ that evolves the polynomial solution within each element over the time step $\Delta t$, hence leading to a one-step time marching scheme. Starting from the known polynomial solution at time $t^n$, the method solves a local Cauchy problem in the small within each element by projecting the governing PDE system onto a space-time basis. This yields a high-order accurate space-time polynomial predictor that approximates the solution in the space-time element $P_i \times \Delta t$. Then, the predictor is extensively used to integrate the governing equations onto a space-time control volume to perform the corrector step, where informations are travelling among control volumes via the computation of numerical fluxes according to the DG formulation. 
%

The local predictor is expressed over each space-time subcell $T_{ik} \times \Delta t$ of $P_{i} \times \Delta t$ as a space-time expansion:
\begin{equation}
	\mathbf{q}_h(\mathbf{x}, t) = \sum_{\ell=1}^{\cT} \theta_{\ell}(\tilde{\xi})\, \hat{\mathbf{q}}_{\ell,ik} = \theta_{\ell}(\tilde{\xi})\, \hat{\mathbf{q}}_{\ell,ik}, 
    \quad \bxit = (\bxi,\tau) = (\xi, \eta, \zeta, \tau),\label{eq:ader00}
\end{equation}
where $\cT$ is the number of space-time degrees of freedom, $\theta_{\ell}(\bxit)$ are the space-time basis functions on the space-time reference element $\hat{T}\times [0,1]$ defined by the space-time coordinates $\bxit = (\xi, \eta, \zeta, \tau) \in \hat{T}\times  [0,1]$, and $\hat{\mathbf{q}}_{\ell,ik}$ are the unknown expansion coefficients. Notice that in the formula above and throughout the entire manuscript we use tensor index notation with Einstein convention which implies summation over repeated indexes. Furthermore, we remark that $\cT$ denotes the number of the space-time degrees of freedom on each sub-tetrahedron, while $\Nst_i$ denotes the total number of space-time degrees of freedom on the entire cell $P_i$ without repetition. The space-time basis functions $\theta_{\ell}$  are defined by means of a tensor product between the finite element basis in space on the reference tetrahedron and the one-dimensional Lagrange interpolation basis functions along the time coordinate passing through the Gauss-Legendre points (see Appendix A in \cite{BOSCHERI2019108899} for explicit formulae of the basis functions on simplex control volumes for two- and three-dimensional problems).

We rewrite the governing PDE in reference coordinates $\bxit$ leading to the transformed system:
\begin{equation}
\frac{\partial \mathbf{Q}}{\partial \tau} + \nabla_{\bxi} \cdot \mathbf{F}^*(\mathbf{Q}, \nabla \mathbf{Q}) = \mathbf{0}, \, 
\label{eq:PDE_ader0}
\end{equation}
\begin{equation*}
\mathbf{F}^*(\mathbf{Q}, \nabla \mathbf{Q}) = (\mathbf{f}^*,\mathbf{g}^*,\mathbf{h}^*),  \qquad 
\nabla_{\bxi} = \left( \frac{\partial}{\partial \xi},\, \frac{\partial}{\partial \eta}, \, \frac{\partial}{\partial \zeta} \right),
\end{equation*}
with the transformed fluxes explicitly given by
\begin{subequations}
\begin{align}
	\mathbf{f}^* &= \Delta t \left( f \frac{\partial \xi}{\partial x} + g \frac{\partial \xi}{\partial y} + h \frac{\partial \xi}{\partial z} \right), 
	\\
	\mathbf{g}^* &= \Delta t \left( f \frac{\partial \eta}{\partial x} + g \frac{\partial \eta}{\partial y} + h \frac{\partial \eta}{\partial z} \right), \\
	\mathbf{h}^* & = \Delta t \left( f \frac{\partial \zeta}{\partial x} + g \frac{\partial \zeta}{\partial y}  + h \frac{\partial \zeta}{\partial z} \right).
\end{align}	
\label{eq:modified_fluxes}
\end{subequations}
Then, the rescaled PDE \eqref{eq:PDE_ader0} is multiplied by a set of test functions $\theta_m$, of the same form of the basis functions $\theta_{\ell}$, and subsequently integrated over the reference space-time element, hence yielding the following weak formulation: 
\begin{equation}
\int\limits_0^1 \int \limits_{\hat{T}} \theta_m(\tilde{\xi}) \left( \frac{\partial \mathbf{q}_h}{\partial \tau} + \nabla_{\bxi} \cdot \mathbf{F}^*(\mathbf{q}_h) \right) \, d\bxit = \mathbf{0}, \quad \forall m \in [ 0, \cT ],\label{eq:ader1}
\end{equation}
which, using the ansatz \eqref{eq:ader00} for the predictor solution, becomes
\begin{equation}
\int\limits_0^1 \int \limits_{\hat{T}} \theta_m \frac{\partial \theta_\ell}{\partial \tau} \hat{\mathbf{q}}_\ell \, d\bxit
+
\int\limits_0^1 \int \limits_{\hat{T}} \theta_m \nabla_{\bxi} \cdot \mathbf{F}^*(\theta_\ell \hat{\mathbf{q}}_\ell) \, d\bxit = \mathbf{0}, \quad \forall m \in [ 0, \cT ]. \label{eq:ader2-1}
\end{equation}
This equation being true for each sub-tetrahedron $T_{ik}$, multiplying by $\abs{J_{ik}}$ the determinant of the spatial Jacobian matrix 
\begin{equation}
    J_{ik}=\frac{\partial \mathbf{x}}{\partial \boldsymbol{\xi}}\Bigg\rvert_{T_{ik}},
\end{equation}
and summing over all the sub-tetrahedra of the cell, we obtain
\begin{equation}
\abs{P_i} \int\limits_0^1 \int \limits_{\hat{T}} \theta_m \frac{\partial \theta_\ell}{\partial \tau} \hat{\mathbf{q}}_{\ell,ik} \, d\bxit =-
\sum_{k=1}^{N_i} \abs{J_{ik}}
\int\limits_0^1 \int \limits_{\hat{T}} \theta_m \nabla_{\bxi} \cdot \mathbf{F}^*(\theta_\ell \hat{\mathbf{q}}_{\ell,ik}) \, d\bxit,  \quad \forall m \in [ 0, \cT ]. \label{eq:ader2}
\end{equation}
The above expression constitutes a nonlinear system for the space-time coefficients $\hat{\mathbf{q}}_\ell$, and it can be compactly written in matrix form as
\begin{equation}
\K_{t,i} \, \hat{\mathbf{q}}_{\ell,i} = -\sum_{k =1}^{N_{i}} \abs{J_{ik}} \left( \K_{\xi} \, \hat{\mathbf{f}}^*_{ik} + \K_{\eta} \, \hat{\mathbf{g}}^*_{ik} + \K_{\zeta} \, \hat{\mathbf{h}}^*_{ik} \right), \label{eq:ader3}
\end{equation}
with the time stiffness matrix
\begin{equation}
    \K_{t,i} = \abs{P_i}\int\limits_0^1 \int \limits_{\hat{T}} \theta_m \frac{\partial \theta_\ell}{\partial \tau} d\bxit,
    \label{eq:ADER-time_stiff}
\end{equation}
and the spatial stiffness matrices
\begin{equation}
\K_{\xi} = \int\limits_0^1 \int \limits_{\hat{T}} \theta_m \frac{\partial \theta_\ell}{\partial \xi} d\bxit, \quad
\K_{\eta} = \int\limits_0^1 \int \limits_{\hat{T}} \theta_m \frac{\partial \theta_\ell}{\partial \eta} d\bxit, \quad
\K_{\zeta} = \int\limits_0^1 \int \limits_{\hat{T}} \theta_m \frac{\partial \theta_\ell}{\partial \zeta} d\bxit.  
\label{eqn.Kspace}
\end{equation}
The numerical fluxes $\hat{\mathbf{f}}^*_{ik}$, $\hat{\mathbf{g}}^*_{ik}$, $\hat{\mathbf{h}}^*_{ik}$ are evaluated exploiting the collocation property of the AFE nodal basis, hence directly from $\hat{\mathbf{q}}_{\ell,i}$ and $\nabla \theta_\ell \hat{\mathbf{q}}_{\ell,i}$. 
Let us underline that the spatial stiffness matrices \eqref{eqn.Kspace} are defined in the reference space-time element $\hat{T}\times[0,1]$, therefore, they do not depend neither on space nor on time, since the space-time dependency is taken into account by the Jacobian $\Delta t \abs{J_{ik}}$. However, the time stiffness matrix $\K_{t,i}$ must handle simultaneously all degrees of freedom $\Nst_i$ for the predictor solution $\hat{\mathbf{q}}_{\ell,i}$, therefore it is polyhedral-dependent and its dimension is necessarily $\Nst_i\times \Nst_i$. Following the development proposed in \cite{BoscheriAFE2022}, an integration by parts in time of the left-hand-side of \eqref{eq:ader3} provides an upwinding approximation in the time direction:
\begin{equation}
    {\K}_{1,i} \hat{\mathbf{q}}_{\ell,i} = \sum_{k =1}^{N_{i}} \abs{J_{ik}} \mathbf{F_0} \hat{\uu}_{\ell,ik}
    - \sum_{k =1}^{N_{i}} \abs{J_{ik}} \left( \K_{\xi} {\hat{\f}}^*_{ik} + \K_{\eta} {\hat{\g}}^*_{ik} + \K_{\zeta} {\hat{\h}}^*_{ik} \right),
    \label{eq:ADER-localCauchy}
\end{equation}
with
\begin{equation*}
{\K}_{1,i} = \abs{P_i} \sum_{k =1}^{N_{i}} \left[
\int \limits_{\hat{T}} \theta_m(\xi,1)\theta_{\ell}(\xi,1) d\bxi
-\int \limits_0^1 \int \limits_{\hat{T}} \frac{\partial \theta_m}{\partial \tau} \theta_{\ell} d\bxit
\right],
\,
\mathbf{F_0} = \int \limits_{\hat{T}} \theta_k(\xi,0)\, \phi_{\ell}(\xi)\, d\bxi.
\end{equation*}
Here, $\mathbf{F_0}$ corresponds to a universal matrix defined on the reference element $\hat{T}$, while $\K_{1,i}$ is an element-wise matrix, as well as the time stiffness matrix $\K_{t,i}$ given in \eqref{eq:ADER-time_stiff}. The element-local nonlinear algebraic system \eqref{eq:ADER-localCauchy} is solved iteratively using a fixed point method for the unknown space-time expansion coefficients $\hat{\q}_{\ell,i}$:
\begin{equation}
\hat{\mathbf{q}}_{\ell,i}^{r+1} = \K_{1,i}^{-1} \left[
\sum_{k =1}^{N_{i}} \abs{J_{ik}} \mathbf{F_0} \hat{\uu}_{\ell,ik}
- \sum_{k =1}^{N_{i}} \abs{J_{ik}} \left( \K_{\xi} {\hat{\f}}^{*,r}_{ik} + \K_{\eta} {\hat{\g}}^{*,r}_{ik} + \K_{\zeta} {\hat{\h}}^{*,r}_{ik} \right)
\right],
\label{eq:ader4}
\end{equation}
with $r$ denoting the iteration number. The iteration process stops when the residual of \eqref{eq:ader4} is lower than a given tolerance (typically set to $10^{-12}$). The convergence of the iterative solver has been proven for linear PDEs in \cite{JACKSON2017409} and for nonlinear PDEs in \cite{Busto2019}.

An efficient way to reduce the computational intensity of the algorithm is to precompute and store the inverse matrix $\K_{1,i}$ for each polyhedral element $P_i$ during the preprocessing stage. Then, all other matrices appearing in the nonlinear equation \eqref{eq:ADER-localCauchy} are evaluated on the reference elements $\hat{T}\times [0,1]$, where the newly introduced agglomerated finite element basis functions are uniquely defined for each subcell $T_{ik}$.

Once the iterative process \eqref{eq:ader4} converged, the time integrated predictor solution $\bar{\mathbf{q}}_h(\mathbf{x})$ as well as its associated physical flux tensor $\overline{\mathbf{F}}_h=(\overline{\mathbf{f}}_h,\overline{\mathbf{g}}_h,\overline{\mathbf{h}}_h)$ can be directly computed:
\begin{subequations}
\begin{align}
\bar{\mathbf{q}}_h(\mathbf{x}) &= \int \limits_{t^n}^{t^{n+1}} \mathbf{q}_h \, dt
= \Delta t \sum_{k =1}^{N_{i}} \mathbf{T}_{\tau} \hat{\mathbf{q}}_{\ell,ik}, \\
\overline{\mathbf{f}}_h(x) &= \int_{t^n}^{t^{n+1}} \hat{\mathbf{f}}_h^{*} \, dt=\Delta t \sum_{k =1}^{N_{i}} \mathbf{T}_{\tau} \hat{\mathbf{f}}^{*}_{\ell,ik}, \\
\overline{\mathbf{g}}_h(x) &= \int_{t^n}^{t^{n+1}} \hat{\mathbf{g}}_h^{*} \, dt=\Delta t \sum_{k =1}^{N_{i}} \mathbf{T}_{\tau} \hat{\mathbf{g}}^{*}_{\ell,ik}, \\
\overline{\mathbf{h}}_h(x) &= \int_{t^n}^{t^{n+1}} \hat{\mathbf{h}}_h^{*} \, dt=\Delta t \sum_{k =1}^{N_{i}} \mathbf{T}_{\tau} \hat{\mathbf{h}}^{*}_{\ell,ik}.
\end{align}
\label{eqn.TAvectors}
\end{subequations}
The universal matrix $ \mathbf{T}_{\tau}$ is evaluated as the integral of the space-time
basis functions over the reference time interval $[0,1]$, that is
\begin{equation}
	\mathbf{T}_{\tau} = \int \limits_0^1 \theta_{\ell}(\tau) \, d\tau.
\end{equation} 
Finally, the time-integrated predictor solution $\bar{\mathbf{q}}_h(\mathbf{x})$ being not time dependent, we can express $\bar{\mathbf{q}}_h(\mathbf{x})$ using the same agglomerated basis functions $\phi_{\ell}(\mathbf{x})$ adopted in \eqref{eq:interp_u} for the numerical interpolation of $\uu_h(\mathbf{x}, t)$:
\begin{equation}
\bar{\mathbf{q}}_h(\mathbf{x}) = \sum_{\ell=1}^{\N_i} \phi_{\ell}(\mathbf{x}) \, \hat{\bar{\mathbf{q}}}_{\ell,i},
\tag{3.27}
\end{equation}
where $\N_i$ corresponds to the total spatial number of degrees of freedom in $P_i$, and $\hat{\bar{\mathbf{q}}}_{\ell,i}$ are the new expansion coefficients of the predictor in the cell.

\subsection{Fully discrete quadrature-free one-step ADER DG scheme}

After computing the time integrated predictor solution $\bar{\mathbf{q}}_h(\mathbf{x})$, a final corrector step is applied to obtain a fully-discrete DG scheme. 
We start from the weak form of the conservation law \eqref{eq:conservation_law} obtained upon multiplication by the test functions $\phi_m$, which are taken of the same form of the basis \eqref{eq:interp_u}, and integration over each space-time element $P_i \times [t^n, t^{n+1}]$:
\begin{equation}
\int \limits_{t^n}^{t^{n+1}} \int \limits_{P_i}
\phi_m \left( \frac{\partial \mathbf{Q}}{\partial t} + \nabla \cdot \mathbf{F}(\mathbf{Q}, \nabla \mathbf{Q}) \right)
\, d\mathbf{x} \, dt = \bzero.
\label{eq:CL-integrated}
\end{equation}
Integrating by parts the flux divergence term as it is generally done for DG schemes leads to
\begin{equation}
\int \limits_{t^n}^{t^{n+1}} \int \limits_{P_i}
\phi_m \frac{\partial \mathbf{Q}}{\partial t}
\, d\mathbf{x} \, dt
+
\int \limits_{t^n}^{t^{n+1}} \int \limits_{\partial P_i}
\phi_m \, \mathbf{F} \cdot \mathbf{n} \, dS \, dt
-
\int \limits_{t^n}^{t^{n+1}} \int \limits_{P_i}
\nabla \phi_m \cdot \mathbf{F} \, d\mathbf{x} \, dt
= \bzero,
\label{eq:CL-integratedIPP}
\end{equation}
where $\partial P_i $ denotes the boundary of the polyhedral cell $P_i$ and $\mathbf{n}$ the corresponding outward pointing normal vector. Using the already time integrated predictor $\bar{\mathbf{q}}_h(\mathbf{x})$ and collocated fluxes computed with \eqref{eqn.TAvectors}, we can then map the integrals to the reference element representation and obtain the quadrature-free one-step and fully-discrete update:
\begin{equation}
\mathbf{M}_i \left( \hat{\mathbf{u}}_{\ell,i}^{n+1} - \hat{\mathbf{u}}_{\ell,i}^{n} \right)
 + \sum_{k=1}^{N_i}
\abs{\partial P_{ik}} \, \mathbf{Z}_{\bxi} \, \mathcal{G}_k\left(\bar{\mathbf{q}}_h^+, \bar{\mathbf{q}}_h^-\right)
=
\sum_{k=1}^{N_i}
\abs{J_{ik}}
\left( \mathbf{V}_\xi \overline{\mathbf{f}}_h + \mathbf{V}_\eta \overline{\mathbf{g}}_h + \mathbf{V}_\zeta \overline{\mathbf{h}}_h \right),
\end{equation}
where $\mathbf{M}_i$ is the element-wise mass matrix
\begin{equation*}
    \mathbf{M}_i = \abs{P_i} \sum_{k=1}^{N_i}
  \int \limits_{\hat{T}} \phi_m \, \phi_\ell \, d\bxi,
\end{equation*}
$\mathbf{V}_\xi$, $\mathbf{V}_\eta$, and $\mathbf{V}_\zeta$ are universal derivative matrices defined on the reference tetrahedron $\hat{T}$, and $\mathbf{Z}_{\bxi}$ is a face mass matrix defined on reference faces $\hat{K} \subset \partial \hat{T}$:
\begin{align*}
  \mathbf{Z}_{\bxi} = \int_{\hat{K}} \phi_m \, \phi_\ell \, d\chi, \quad
  \mathbf{V}_\xi = \int \limits_{\hat{T}} \frac{\partial \phi_m}{\partial \xi} \, \phi_\ell \, d\bxi, \\
  \mathbf{V}_\eta = \int \limits_{\hat{T}} \frac{\partial \phi_m}{\partial \eta} \, \phi_\ell \, d\bxi, \quad
  \mathbf{V}_\zeta = \int \limits_{\hat{T}} \frac{\partial \phi_m}{\partial \zeta} \, \phi_\ell \, d\bxi.
\end{align*}
At each face $\partial P_{ik}$ of surface $\abs{\partial P_{ik}}$, the numerical flux $\mathcal{G}_k$ is computed using the right and left states $\left(\bar{\mathbf{q}}_h^+, \bar{\mathbf{q}}_h^-\right)$ which are available from the predictor step. The Rusanov--type numerical flux function is adopted, with its extension to viscous terms according to \cite{Gassner2007}, hence yielding
\begin{equation}
\mathcal{G}_k = \frac{1}{2} \left( \overline{\mathbf{F}}_h^+ + \overline{\mathbf{F}}_h^- \right) \cdot \mathbf{n}
- \frac{1}{2} \left( \abs{\boldsymbol{\lambda}_{\max}} + 2\,\varepsilon\, \abs{\boldsymbol{\lambda}_{\max}^{v}} \right) \left( \overline{\q}_h^+ - \overline{\q}_h^- \right),
\tag{3.33}
\end{equation}
where $$\varepsilon = \frac{2N + 1}{(h^+ + h^-)\sqrt{\frac{\pi}{2}}} $$ is an artificial viscosity coefficient, with $h^+$ and $h^-$ being characteristic lengths of the element $P_i$ and its neighbor $P_j$. The quantities $\abs{\boldsymbol{\lambda}_{\max}}$ and $\abs{\boldsymbol{\lambda}_{\max}^{v}}$ are computed by taking the maximum absolute values of the convective and viscous eigenvalues from \eqref{eq:eigenvalues3D} between the states $\overline{\q}_h^+$ and $\overline{\q}_h^-$, and then projecting them along the face normal $\mathbf{n}$. 

\subsection{Artificial viscosity limiter}

High-order discontinuous Galerkin (DG) methods are prone to non-physical oscillations near discontinuities such as shock waves due to the Gibbs phenomenon. To maintain stability in the presence of shocks or steep gradients, we must incorporate a limiter technique. Among the various approaches proposed in the literature \cite{Boscheri2016, Boscheri2019, BOSCHERI2019108899, Rosa2017, DUMBSER201447, FRANK2020112665, KUZMIN20103077, Markert2021, Sonntag2014, ZANOTTI2015204, ZHU2020109105} (ranging from classical a priori limiters to modern a posteriori techniques, including bound-preserving strategies and multi-level refinement approaches), we adopt here a simple artificial viscosity method \cite{BASSI2018186,GASSNER2015247, Klockner2011, Massa2020, Persson2006}, applied exclusively to elements detected as \textit{“troubled”}. Consequently, the limiter procedure consists of two stages: (i) detection, and (ii) limiting. 
This approach ensures that the limiter selectively adds dissipation only where needed (e.g., near shocks or strong gradients), while preserving the high-order accuracy of the DG scheme in smooth regions. In three dimensions, the method naturally accounts for contributions from all velocity components and for the additional geometrical complexity of polyhedral faces, providing a robust stabilization strategy for fully unstructured three-dimensional meshes.

To identify troubled cells in three-dimensional polyhedral meshes, we employ a flattener variable $\beta_i$, following the approach of \cite{BALSARA20127504}, computed at the beginning of each time step. The detector compares the divergence of the velocity field $\nabla \cdot \mathbf{v}$ within each element $P_i$ to the minimum sound speed $c = \sqrt{\gamma R T}$ among $P_i$ itself and its direct neighbors. The discrete divergence and the minimum sound speed are evaluated as
\begin{equation}
\nabla \cdot \mathbf{v} \big\rvert_i = \frac{1}{\abs{P_{i}}} \sum_{k=1}^{N_i} \abs{\partial P_{ik}}  (\mathbf{v}^+ - \mathbf{v}^-) \cdot \mathbf{n}, 
\qquad
c_{s,\min} = \min_{1\leq k \leq N_i} (c_s^+, c_s^-),
\label{eq:divergence3D}
\end{equation}
where $(+,-)$ indicate the right and left states across face $\partial P_{ik}$ of $P_i, \, \forall k \in \{ 1,N_i \}$. The flattener variable is then defined as
\begin{equation}
\beta_i = \min \Big[ 1, \max \Big( 0, - \frac{\nabla \cdot \mathbf{v}\rvert_i + m_1 c_{s,\min}}{m_1 c_{s,\min}} \Big) \Big],
\label{eq:beta3D}
\end{equation}
with $m_1 = 0.1$ as proposed in \cite{BALSARA20127504}. If the detector is activated on cell $P_i$, the element is therefore flagged as troubled and an artificial viscosity $\mu_{\text{add},i}$ is added to the physical viscosity $\mu$, leading to a local (or cell-dependent) viscosity coefficient $\mu_i = \mu + \mu_{\text{add},i}$ used in the viscous stress tensor \eqref{eq:stress_tensor}. The artificial viscosity coefficient is determined such that a unity mesh Reynolds number $Re_i$ is assigned to the troubled cells, that is
\begin{equation}
Re_i = \frac{\rho_i \, \abs{\boldsymbol{\lambda}_{\max}}_i \, h_s}{\mu_i}, 
\qquad 
h_s = \frac{h_i}{2N+1},
\end{equation}

\section{Numerical tests} \label{sec:num_res}

In this section, we assess the performance of the proposed high-order ADER-DG scheme on three-dimensional unstructured polyhedral meshes by running a comprehensive set of numerical experiments. The selected test cases have been carefully designed to evaluate various aspects of the method, including accuracy on smooth flows through convergence tests, robustness in the presence of shocks, and the ability to handle viscous diffusion. Whenever possible, numerical results are compared against exact analytical solutions or well-established reference data (\cite{BOSCHERI2022127457,BoscheriAFE2022} and references therein) to demonstrate the accuracy, stability, and robustness of the proposed ADER-DG method for complex three-dimensional flows.

We set $\text{CFL}=0.3$ in all simulations, and the ratio of specific heats is $\gamma=1.4$, the gas constant is $R=1$,
thus the specific heat capacity at constant volume is $c_v = 2.5$. Heat conduction is neglected, hence setting $\kappa=0$.


\subsection{Two-dimensional Steady Vortex} 
The steady vortex test case \cite{HuShuTri} is a classical benchmark to evaluate the spatial accuracy and stability of the numerical scheme under stationary flow conditions. The vortex configuration remains unchanged in time, allowing to isolate and quantify spatial discretization errors without temporal contamination. This test is particularly useful to verify that the proposed method can maintain high-order accuracy on complex unstructured meshes. We consider the domain $\Omega = [0, 10]^3$, discretized using unstructured polyhedral meshes. The two-dimensional vortex center is initialized at $(5,5,z),z\in[0,10]$, and the initial conditions correspond to a smooth and divergence-free vortex flow given by
\begin{subequations}
	\begin{align}
\rho(\mathbf{x}, 0) =& 1+\delta\rho, \\
\mathbf{u}(\mathbf{x}, 0) =& (\delta u, \delta v, 0),  \\
p(\mathbf{x}, 0) =& 1+\delta p,
  \end{align}
\label{eq:vortex1}
\end{subequations}
where the perturbations explicitly read
\begin{subequations}
\begin{align}
    \delta p(\mathbf{x},0) =& (1+\delta \cT)^{\frac{\gamma}{\gamma-1}}-1, \\ \, 
    \delta \rho(\mathbf{x},0) =& (1+\delta \cT)^{\frac{1}{\gamma-1}}-1, \\ \, 
    \begin{pmatrix}
    \delta u \\
    \delta v \\
    \delta w
    \end{pmatrix}(\mathbf{x},0) =& \frac{\epsilon}{2\pi}e^{\frac{1-\gamma^2}{2}}
    \begin{pmatrix}
    -(y-5) \\
    \phantom{+}(x-5) \\
    0
    \end{pmatrix}.
\end{align}
\label{eq:vortex2}
\end{subequations}
The vortex strength is $\epsilon=5$, and $\delta \cT$ denotes the temperature perturbation
\begin{align*}
    \delta \cT(x,y,z,t=0) =& - \frac{(\gamma-1)\epsilon^2}{8\gamma\pi^2}e^{1-r^2},
\end{align*}
with the radial coordinate $r(x,y) = \sqrt{(x-5)^2+(y-5)^2}$.
We impose Dirichlet boundary conditions on all sides. The simulation is run until the final time $t_f = 1$ and the analytical solution is given by the initial condition.

\begin{figure}[!htbp]
    \centering
    \includegraphics[trim=0cm 0cm 0cm 2.32cm, clip, height=0.4\textwidth]{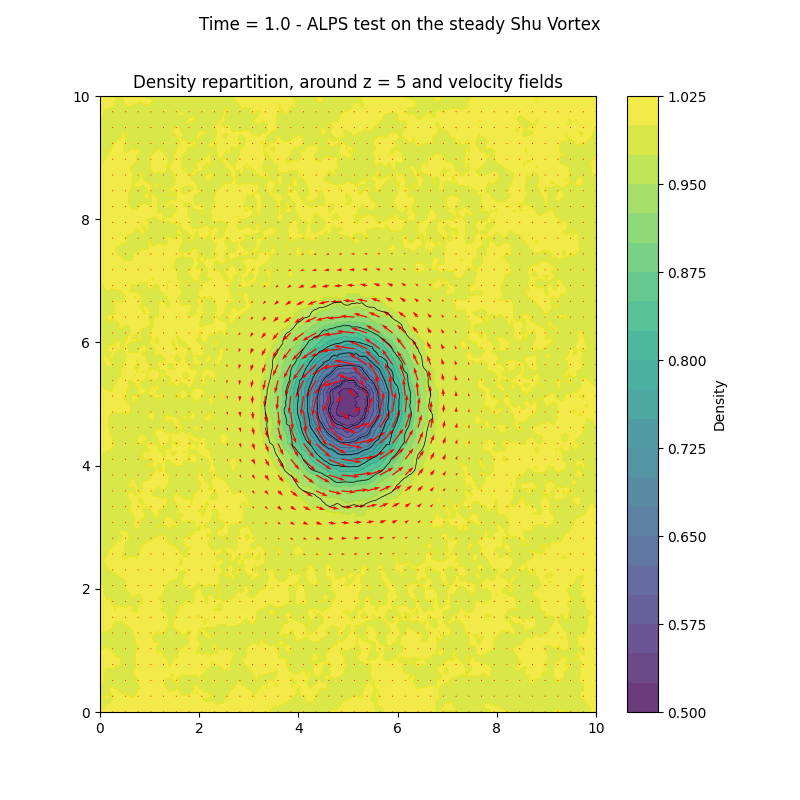}
    \includegraphics[trim=0cm 0cm 0cm 1.74cm, clip, height=0.4\textwidth]{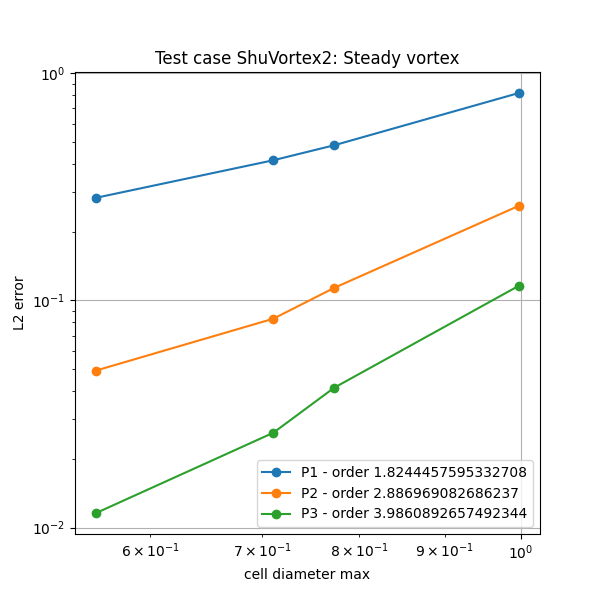}
    \caption{Left: {Density and velocity field of the stationary vortex at $t_f = 1$}, $xy$ plane view at $z=5$. Right: Convergence rates for different polynomial degrees.}
    \label{fig:steady_vortex}
\end{figure}

To assess convergence, we consider a sequence of successively refined unstructured meshes, maintaining a similar element quality. Simulations are performed for different polynomial degrees $N = \{1, 2, 3\}$ of the novel ADER-AFE-DG schemes. The numerical solution is compared to the exact one using $L_1$ and $L_2$ norms:
\begin{equation}
\norm{e}_{L_p} = \left( \sum_{K \in \mathcal{T}_h} \int_K \abs{\rho_h - \rho_e}^p \, d\mathbf{x} \right)^{1/p}, \quad p = \{1, 2\}.
\end{equation}
Table \ref{tab:vortex_errors} reports the errors and the corresponding observed convergence rates, which achieve the expected order $O(h^{N+1})$. Fig. \ref{fig:steady_vortex} shows a slice of the numerical solution at the final time of the simulation and a plot of the convergence rates obtained for three different polynomial degrees.
\begin{table}[!htbp]
    \centering
    \begin{tabular}{|c || c |c || c |c || c |c |}
        \hline
        $h_{\text{max}}$ & $\P_1$ $\norm{e}_{L_1}$ & $\P_1$ $\norm{e}_{L_2}$ & 
        $\P_2$ $\norm{e}_{L_1}$ & $\P_2$ $\norm{e}_{L_2}$ & 
        $\P_3$ $\norm{e}_{L_1}$ & $\P_3$ $\norm{e}_{L_2}$ \\
        \hline \hline
        0.9971   &   16.2025   &   0.8189   &   5.2472   &   0.261   &   1.7093   &   0.116 \\ 
        0.7729   &   8.9488   &   0.482   &   2.1213   &   0.1135   &   0.5394   &   0.0412 \\ 
        0.7108   &   6.8654   &   0.4141   &   1.4248   &   0.0831   &   0.3129   &   0.0261 \\ 
        0.5567   &   4.6514   &   0.2831   &   0.7473   &   0.0491   &   0.1314   &   0.0116 \\ \hline \hline 
        Order   &   2.063   &   1.8451   &   3.1252   &   2.7292   &   4.4486   &   3.9861 \\
        \hline
    \end{tabular}
    \caption{Error norms and convergence rates for the steady vortex test at $t_f = 1$.}
    \label{tab:vortex_errors}
\end{table}

\subsection{Two-dimensional Travelling Vortex}

The travelling vortex problem extends the previous test by introducing a background advection field of velocity $\mathbf{u}=(1,1,0)$. Therefore, space-time accuracy can be empirically checked and the reference solution is given by the time-shifted initial condition \eqref{eq:vortex1}-\eqref{eq:vortex2}. 

\begin{figure}[!htbp]
    \centering
    \includegraphics[trim=0cm 0cm 0cm 2.32cm, clip, height=0.4\textwidth]{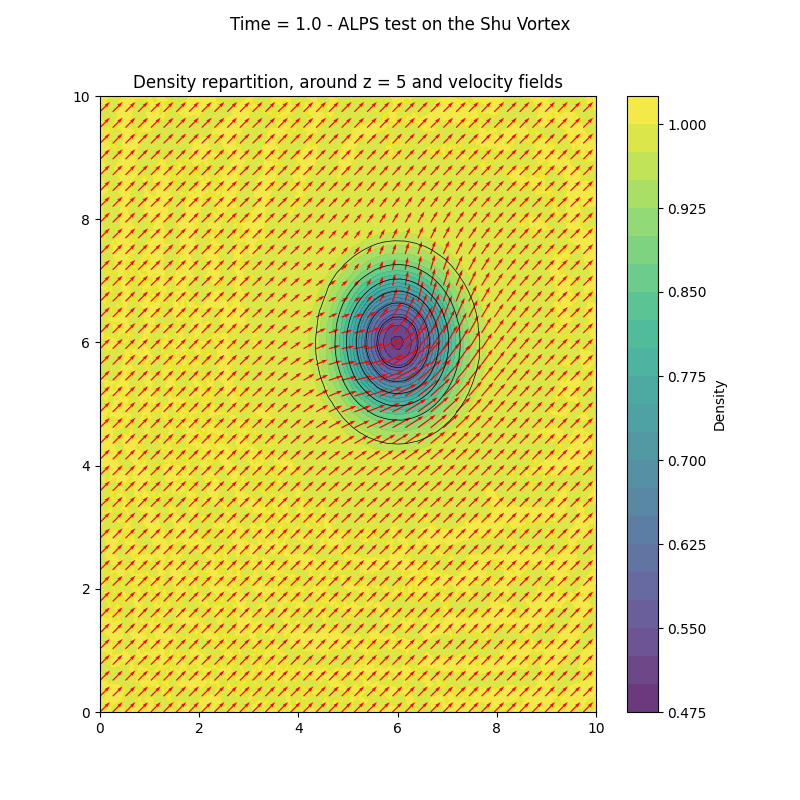}
    \includegraphics[trim=0cm 0cm 0cm 1.74cm, clip, height=0.4\textwidth]{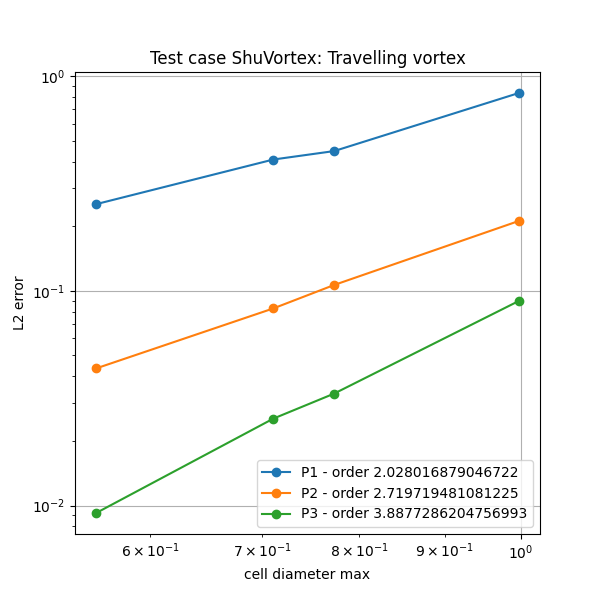}
    \caption{Left: {Density and velocity field of the travelling vortex at $t_f = 1$}, $xy$ plane view at $z=5$. Right: Convergence rates for different polynomial degrees.}
    \label{fig:travelling_vortex}
\end{figure}

%
Fig. \ref{fig:travelling_vortex} shows the density and velocity fields at final time $t_f=1$. Table \ref{tab:travelling_vortex_errors} reports the $L_1$ and $L_2$ errors, as well as convergence rates. 
Again, the method achieves the expected order of convergence $O(h^{N+1})$ in space-time.
\begin{table}[!htbp]
    \centering
    \begin{tabular}{|c || c |c || c |c || c |c |}
        \hline
        $h_{\text{max}}$ & $\P_1$ $\norm{e}_{L_1}$ & $\P_1$ $\norm{e}_{L_2}$ & 
        $\P_2$ $\norm{e}_{L_1}$ & $\P_2$ $\norm{e}_{L_2}$ & 
        $\P_3$ $\norm{e}_{L_1}$ & $\P_3$ $\norm{e}_{L_2}$ \\
        \hline \hline
        0.9971   &   17.6128   &   0.8347   &   4.5751   &   0.2116   &   1.6962   &   0.0897 \\
        0.7729   &   8.6545   &   0.4477   &   2.1719   &   0.1065   &   0.5822   &   0.0332 \\
        0.7108   &   7.1485   &   0.4089   &   1.5766   &   0.0828   &   0.41   &   0.0254 \\
        0.5567   &   4.4258   &   0.2534   &   0.8247   &   0.0435   &   0.152   &   0.0092 \\ \hline \hline
        Order   &   2.3715   &   2.028   &   2.9591   &   2.7197   &   4.1407   &   3.8877  \\
        \hline
    \end{tabular}
    \caption{Error norms and convergence rates for the travelling vortex test at $t_f = 1$.}
    \label{tab:travelling_vortex_errors}
\end{table}

\subsection{First Problem of Stokes}
\begin{figure}[!htbp]
    \centering
    \includegraphics[trim=1.7cm 0.5cm 2.5cm 1.7cm, clip, width=0.49\textwidth]{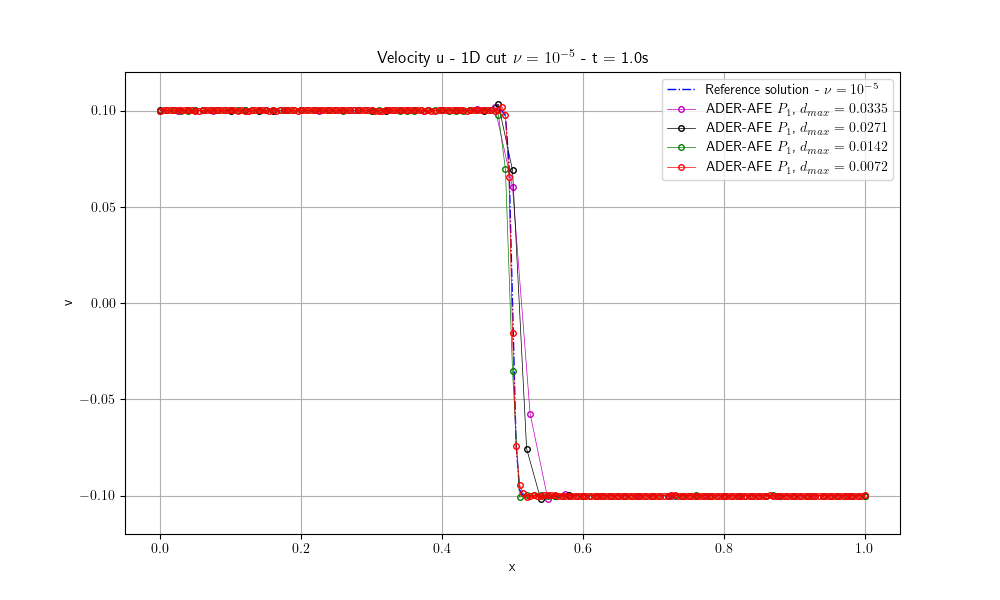}
    \includegraphics[trim=1.7cm 0.5cm 2.5cm 1.7cm, clip, width=0.49\textwidth]{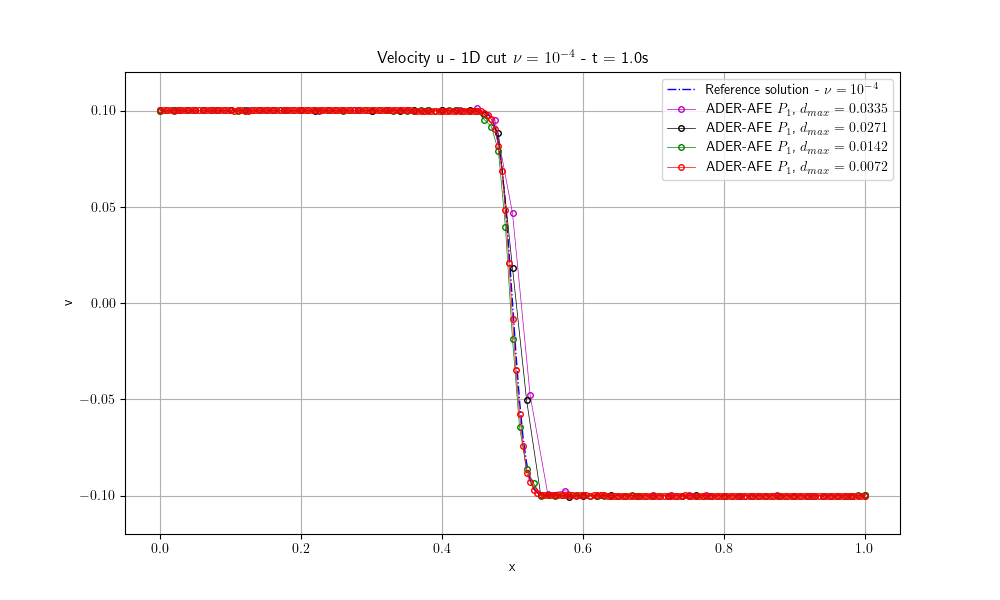}
    \includegraphics[trim=1.7cm 0.5cm 2.5cm 1.7cm, clip, width=0.5\textwidth]{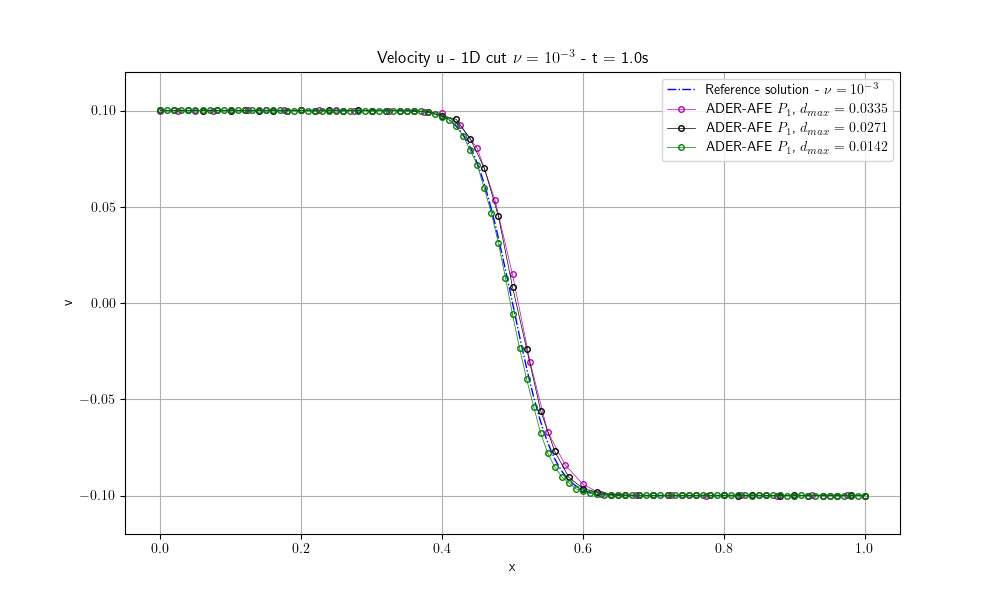}
    \caption{Comparison between exact and numerical solution for the first problem of Stokes at $t_f = 1$ along 
    the $\vec{x}-$direction (1D cut at $y = z = 0.1$), with $\P_1$ elements and different mesh sizes. 
    Viscosity $\nu \in \{ 10^{-5},10^{-4}, 10^{-3}\}$ at respectively top-left, top-right and bottom panels.}
    \label{fig:stokes_profile-P1}
\end{figure}

\begin{figure}[!htbp] 
    \centering
    \includegraphics[trim=1.7cm 0.5cm 2.5cm 1.7cm, clip, width=0.49\textwidth]{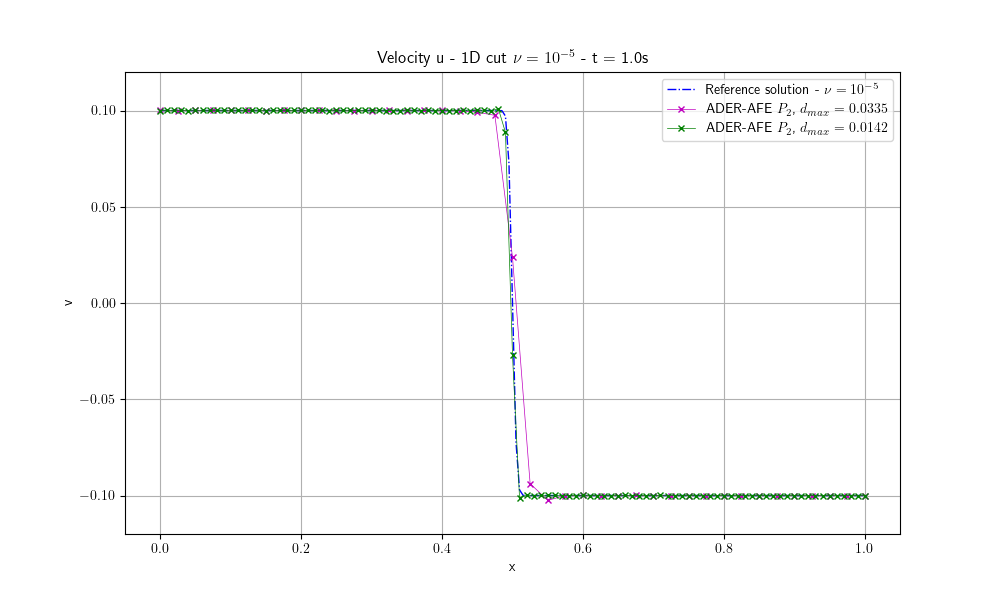}
    \includegraphics[trim=1.7cm 0.5cm 2.5cm 1.7cm, clip, width=0.49\textwidth]{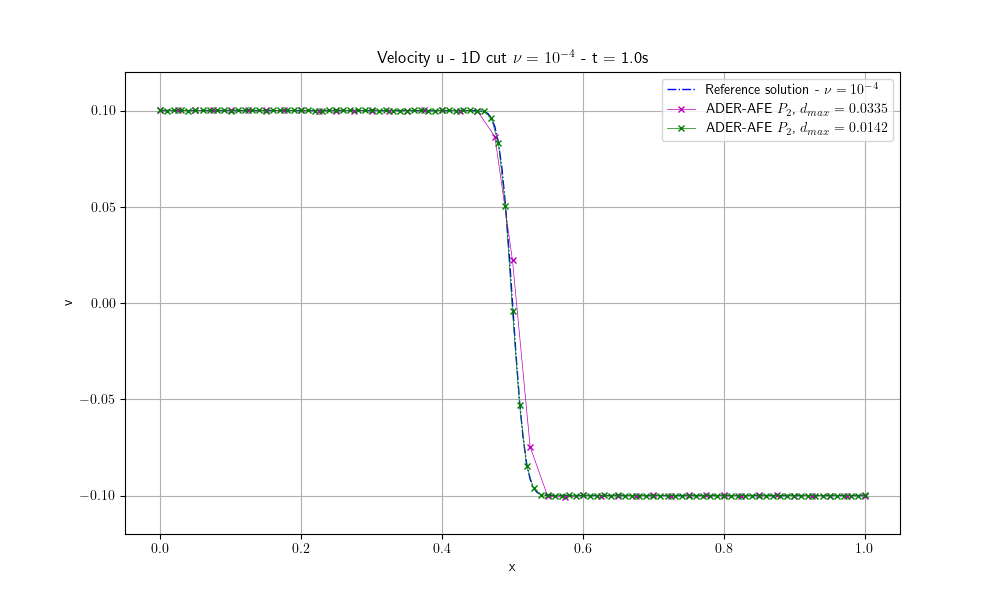}
    \caption{Comparison between exact and numerical solution for the first problem of Stokes at $t_f = 1$ along 
    the $x-$direction (1D cut at $y = z = 0.1$), with $\P_2$ elements and different mesh sizes. 
    Viscosity $\nu \in \{ 10^{-5},10^{-4}\}$ at respectively left and right.}
    \label{fig:stokes_profile-P2}
\end{figure}

This test case \cite{schlichting1997grenzschicht} is designed to validate the correct implementation of viscous terms. By focusing on a purely diffusive flow, it highlights the scheme’s ability to resolve smooth gradients and dissipative phenomena. We consider an infinite incompressible shear layer in computational domain $\Omega = [0,1]\times[0,0.2]^2$. The initial state is defined as
\begin{equation*}
\rho(\mathbf{x}, 0) = 1, \quad u(\mathbf{x}, 0) = w(\mathbf{x}, 0) = 0, \quad p(\mathbf{x}, 0) = \frac{1}{\gamma},
\end{equation*}
and
\begin{equation*}
v(x,0) =
\begin{cases}
v_L = \phantom{+}0.1, & \text{if } x < 0.5, \\
v_R = -0.1, & \text{otherwise}.
\end{cases}
\end{equation*}
For the first problem of Stokes, the analytical solution is given by
\begin{equation*}
v(x,t) = \frac{v_L + v_R}{2} + \frac{v_R - v_L}{2}  \mathrm{erf}\left( \frac{x - 0.5}{2\sqrt{\nu (t+0.05)}} \right),
\end{equation*}
with $\nu$ being the kinematic viscosity. We impose Dirichlet boundary conditions in the $x-$direction, while we impose the exact solution on the remaining sides. The flow is initialized with the above profile and evolved up to $t_f = 1$ with a set of different viscosity parameters $\nu \in \{ 10^{-5},10^{-4}, 10^{-3}\}$. \\
Fig. \ref{fig:stokes_profile-P1} shows the comparison between numerical and analytical solutions at the final time. We can observe the good accuracy of the numerical solutions and their convergence toward the reference profile as the mesh gets finer.
To confirm the good implementation of viscous terms for higher order, Fig. \ref{fig:stokes_profile-P2} presents numerical results at $t_f = 1$ for $\nu \in \{ 10^{-5},10^{-4}\}$ using third-order elements with two different meshes. Once again, the numerical solutions demonstrate good accuracy converging towards the reference solution as the computational mesh gets refined.

\subsection{Three-dimensional Spherical Explosion Problem}

\begin{figure}[!htbp]
    \centering
    \includegraphics[height=0.42\textwidth]{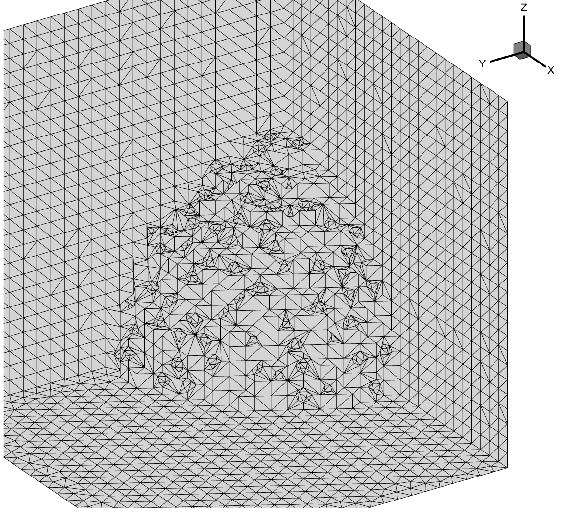}   
    \includegraphics[height=0.42\textwidth]{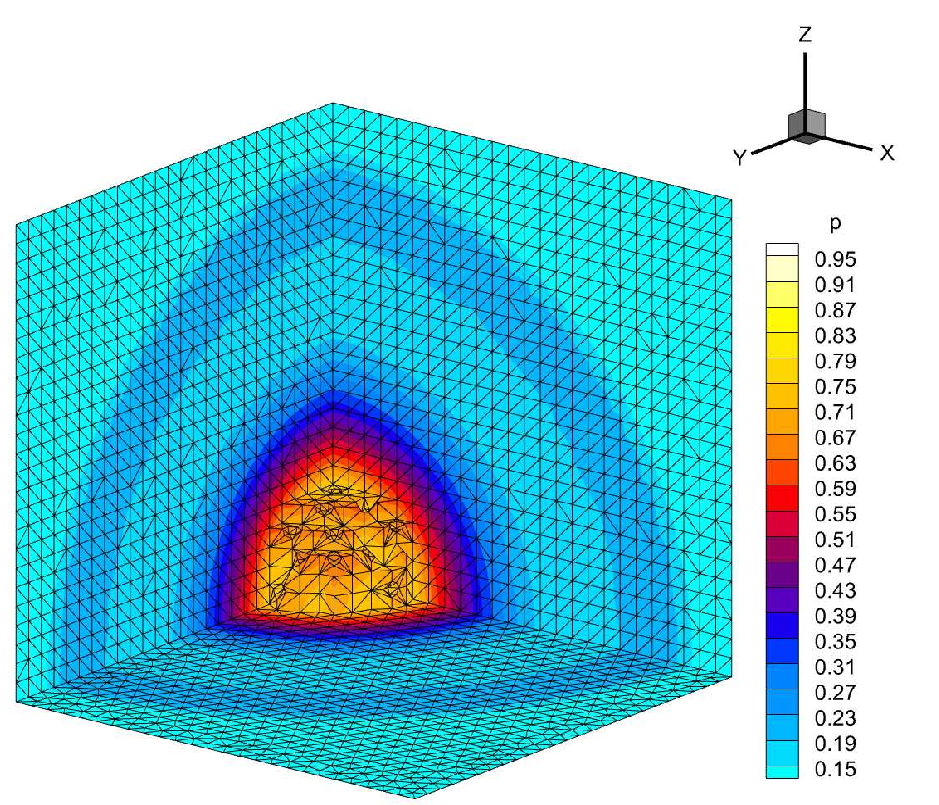}   
    \caption{Left: zoom on the AFE mesh, composed of tetrahedral elements, resulting from the subdivision of polygonal cells. Right: three-dimensional view of the pressure at $t_f = 0.25$ - mesh: $h_{\text{max}}=0.0144$. For both figures, a cut is done at $r=0.2$, which corresponds approximately to the inner boundary of the rarefaction fan at the final time.} 
    \label{fig:sod3d-3d} 
\end{figure}

\begin{figure}[!htbp]
    \centering
    \includegraphics[trim=4.5cm 28.5cm 4.5cm 5.35cm, clip, width=0.99\textwidth]{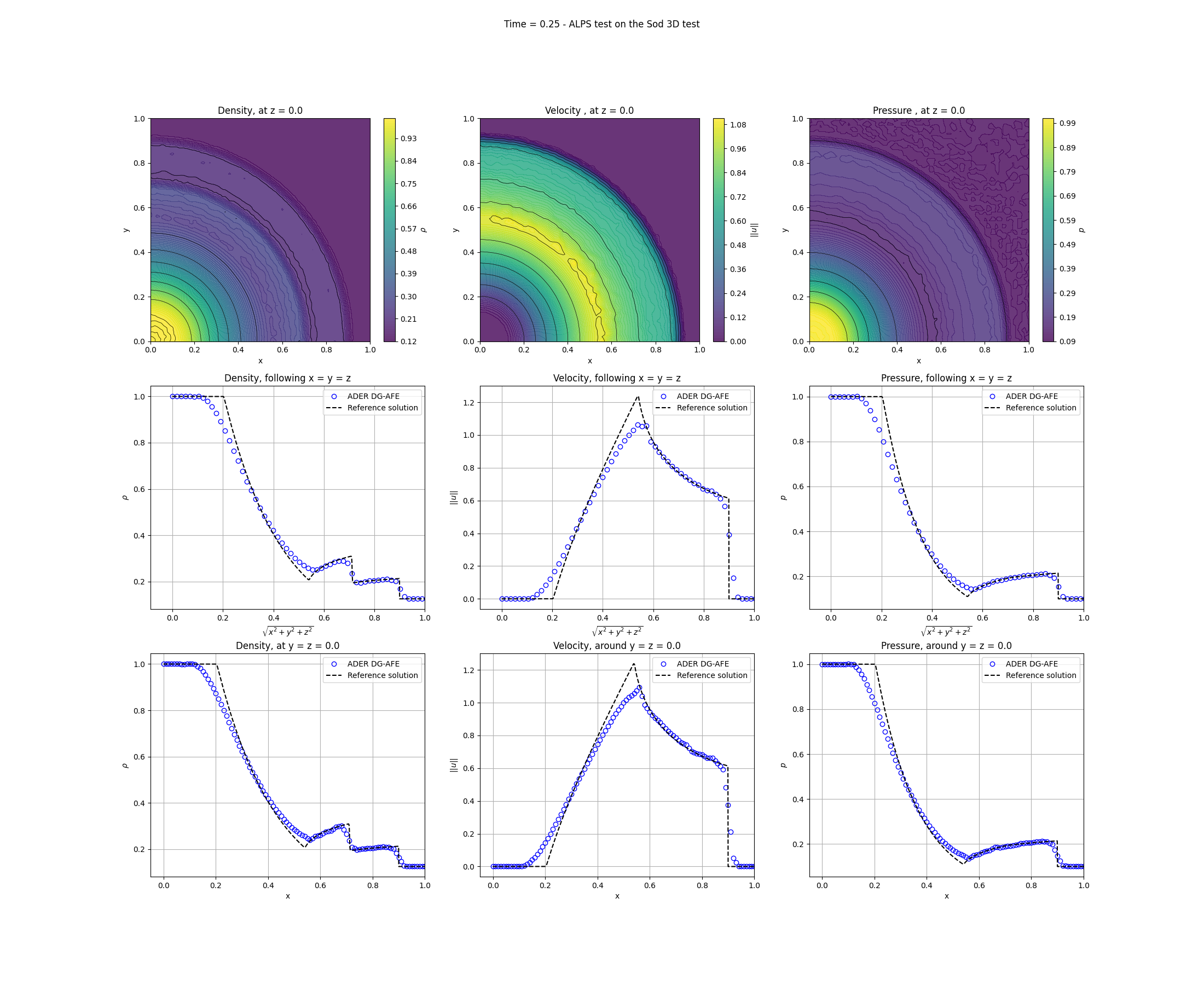}   
    \caption{Two-dimensional contour plots of density, velocity magnitude, and pressure in the plane $z=0$ at $t_f = 0.25$ - mesh: $h_{\text{max}}=0.0287$.} 
    \label{fig:sod3d-2d}
\end{figure}

\begin{figure}[!htbp]
    \centering
    \includegraphics[trim=4.5cm 16cm 4.5cm 17.8cm, clip, width=0.99\textwidth]{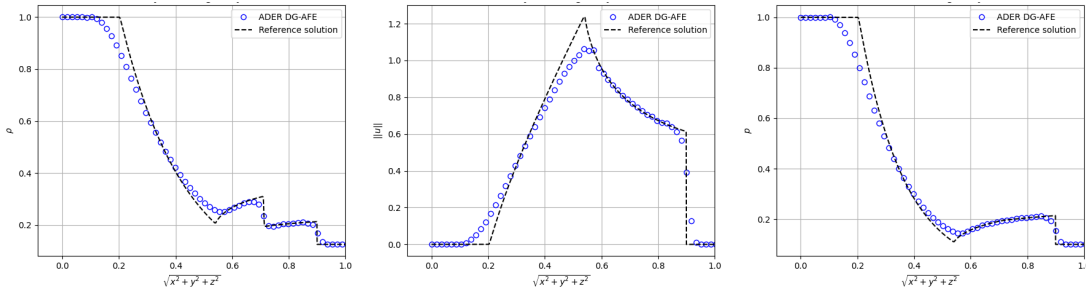}   
    \caption{One-dimensional cuts along the radial direction ($x=y=z$) of density, velocity magnitude, and pressure in the plane $z=0$ at $t_f = 0.25s$ - mesh: $h_{\text{max}}=0.0287$.} 
    \label{fig:sod3d-1d}
\end{figure}

The spherical explosion problem is a very interesting test case involving one shock wave propagating radially from an initial discontinuity as well as a contact and a rarefaction wave. This scenario evaluates the robustness and shock-capturing ability of the numerical method in three dimensions. The artificial viscosity and the shock sensor are therefore activated to stabilize the solution near discontinuities. 
The domain is the cube $\Omega = [0,1]^3$, and the initial condition is given by
\begin{equation*}
\rho(\mathbf{x}, 0) = 
    \begin{cases}
        \rho_{\text{in}}, & \text{if } r < 0.5 \\
        \rho_{\text{out}}, & \text{otherwise}
    \end{cases}, 
\quad \mathbf{u} = 0,
\quad p(\mathbf{x}, 0) = 
    \begin{cases}
        p_{\text{in}}, & \text{if } r < 0.5 \\
        p_{\text{out}}, & \text{otherwise}
    \end{cases},
\end{equation*}
with $r = \norm{\mathbf{x}-\mathbf{x}_C}$, $\mathbf{x}_C=(0.5,0.5,0.5)$. Parameters are set up as $(p_{\text{in}}, p_{\text{out}}) = (1,0.1)$, $(\rho_{\text{in}},\rho_{\text{out}}) = (1,0.125)$. As in \cite{BoscheriAFE2022}, in order to avoid non-physical oscillations at the initial time, the initial condition is slightly smoothed as follows
\begin{equation}
    P(\mathbf{x},0) = (P_{\text{out}} + P_{\text{in}}) + (P_{\text{out}} - P_{\text{in}}) \, \mathrm{erf} \left( \frac{r - R}{\alpha_0} \right), 
    \quad \alpha_0 = 1.5 \times h_{\text{min}},
    \label{eq:initial_smoothing}
\end{equation}
with $h_{\text{min}}$ being the minimal mesh (polyhedral) size, and $P=(\rho,\mathbf{u} , p)$. Thanks to the spherical symmetry of the problem, the solution can be compared with an equivalent one-dimensional problem in the radial direction $r$ (see \cite{Toro2009}). The reference solution at the final time $t_f = 0.25$ is obtained by solving the compressible Euler equations using a second-order MUSCL scheme with the Rusanov flux on a one-dimensional mesh of $15,000$ points over the radial interval $r \in [0,1]$ as proposed in \cite{BoscheriAFE2022}. 
As a result, a spherical shock wave forms and propagates outward, leaving a rarefied low-density region behind. Due to symmetry, the solution remains radially symmetric if the mesh resolution allows it. Fig. \ref{fig:sod3d-3d} illustrates the three-dimensional state of the explosion at the final time: we can see the AFE mesh, composed of tetrahedral elements, resulting from the subdivision of polyhedral cells, and the color map representing the pressure field. A cut is applied at $r=0.2$, which corresponds approximately to the inner boundary of the rarefaction fan at the final time.
Furthermore, Fig. \ref{fig:sod3d-2d} and Fig. \ref{fig:sod3d-1d} depict the numerical solution at the final time through respectively two- and one-dimensional cuts obtained on a finer mesh. In particular, Fig. \ref{fig:sod3d-2d} shows two-dimensional contour plots of density, velocity magnitude, and pressure in the $(x,y)$ plane at $z=0$. We can observe an excellent preservation of radial symmetry. The rarefaction wave, contact discontinuity, and outer shock are clearly identifiable and located at the expected position compared to the reference solution. No instability is observed, nevertheless, a slight granular pattern appears near the outer shock, especially in the pressure field, due to limited mesh resolution in regions with strong gradients. The results shown in Fig. \ref{fig:sod3d-1d} are also in very good agreement with the reference solution. The ADER-AFE-DG scheme shows a controlled dissipation without spurious oscillations.

\subsection{Taylor--Green Vortex}

\begin{figure}[!htbp]
    \centering
    \includegraphics[trim=28.5cm 3cm 7.5cm 4cm, clip, width=0.99\textwidth]{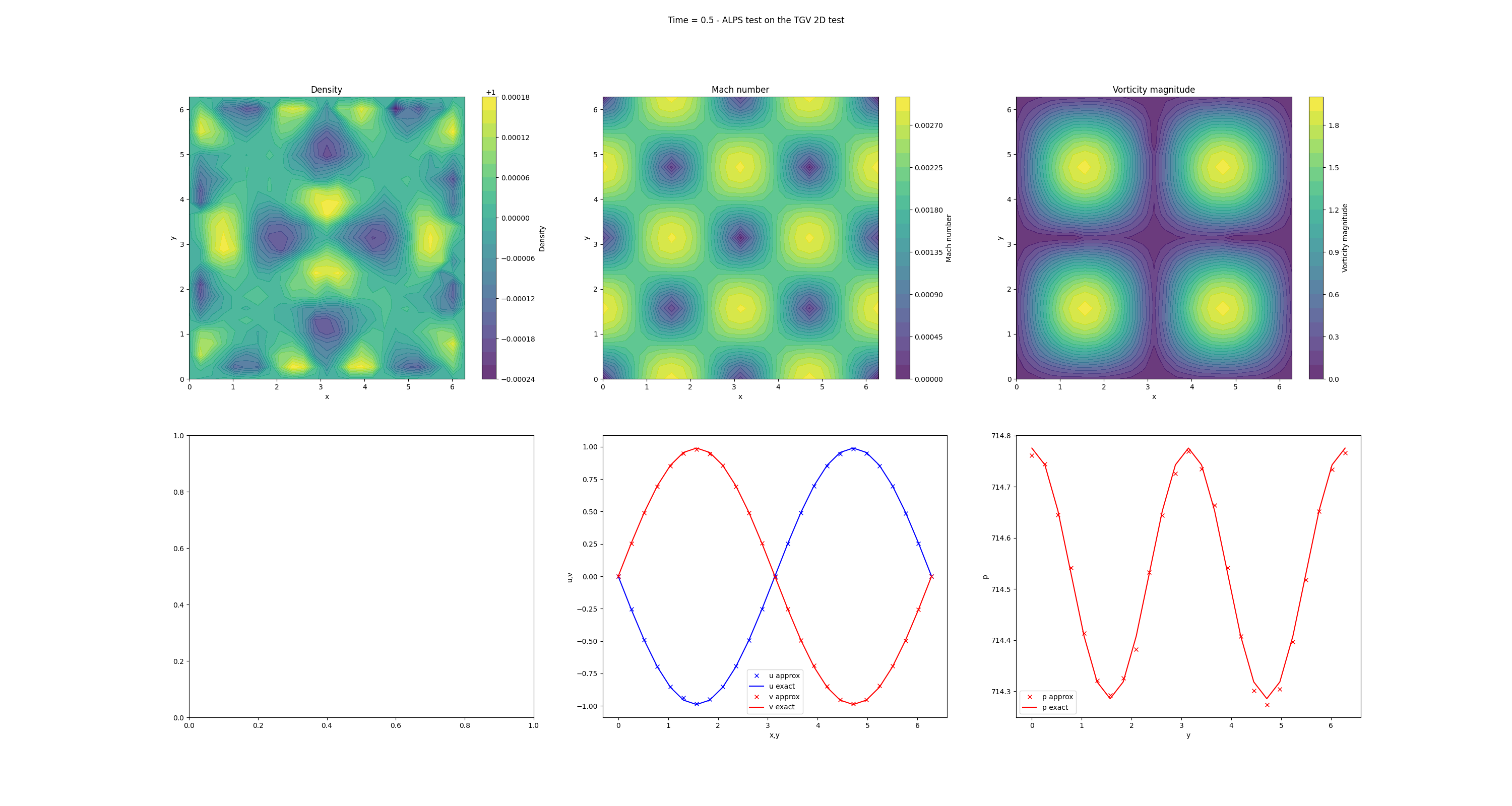}
    \caption{Two-dimensional Taylor-Green vortex problem at $t_f = 0.5$ with $\nu = 10^{-2}$ (and $h_{\text{max}}=0.1804$). Top: Mach number (left) and vorticity magnitude (right) distributions at $z=\pi$. Bottom: One-dimensional cuts along $x-$ and $y-$axis for velocity components $u$ and $v$ as well as for the pressure at respectively $y=\pi$ and $x=\pi$.}
    \label{fig:tgv2d}
\end{figure}

The two-dimensional Taylor-Green vortex is a classical unsteady flow with an exact analytical solution, used to test the convergence and accuracy of numerical schemes for viscous incompressible flows. We consider the two-dimensional vortex in a three-dimensional cubic domain $\Omega = [0, 2\pi]^3$ with Dirichlet boundary conditions, and the initial condition \cite{Taylor1937MechanismOT} writes
\begin{align*}
\rho(\mathbf{x},0) &= 1, \\
u(\mathbf{x},0) &= +\sin(x) \cos(y)e^{-2\nu t},\\
v(\mathbf{x},0) &= -\cos(x) \sin(y)e^{-2\nu t}, \\
w(\mathbf{x},0) &= 0, \\
p(\mathbf{x},0)&=\frac{p_0}{\gamma} + \frac{1}{4} \left( \cos(2x) + \cos(2y) \right)e^{-4\nu t},
\end{align*}
where $\gamma = 1.4$ and $\nu$ is the kinematic viscosity. Fig. \ref{fig:tgv2d} shows the numerical solution for density, Mach number and vorticity magnitude obtained at the final time $t_f =0.5$.
%
The $L_1$ and $L_2$ errors are reported in Table \ref{tab:tgv2d_10-2} and Table \ref{tab:tgv2d_10-5} for $\nu=10^{-2}$ and $\nu=10^{-5} $ respectively. They confirm the correct resolution of the diffusive process and highlights the improvement in numerical accuracy when comparing the errors obtained with $\P_1$ and $\P_2$ elements. The expected order of accuracy is observed.

\begin{table}[htbp]
    \centering
    \begin{tabular}{|c || c |c || c |c |}
        \hline
        $h_{\text{max}}$ & $\P_1$ $\norm{e_u}_{L_1}$ & $\P_1$ $\norm{e_u}_{L_2}$ & 
        $\P_2$ $\norm{e_u}_{L_1}$ & $\P_2$ $\norm{e_u}_{L_2}$  \\
        \hline \hline
        0.4856   &   12.2968   &   0.6186   &   1.2849   &   0.0673 \\ 
        0.4466   &   8.6479   &   0.4383   &   0.8799   &   0.0471 \\ 
        0.3498   &   5.5449   &   0.2914   &   0.4206   &   0.0231 \\ \hline \hline 
        Order   &   2.0591   &   2.0253   &   3.3182   &   3.1845 \\ \hline 
    \end{tabular}\vspace*{0.2cm}
    \begin{tabular}{|c || c |c || c |c |}
        \hline
        $h_{\text{max}}$ & $\P_1$ $\norm{e_v}_{L_1}$ & $\P_1$ $\norm{e_v}_{L_2}$ & 
        $\P_2$ $\norm{e_v}_{L_1}$ & $\P_2$ $\norm{e_v}_{L_2}$   \\
        \hline \hline
        0.4856   &   11.8322   &   0.5984   &   1.2406   &   0.0645 \\
        0.4466   &   8.5964   &   0.4445   &   0.9047   &   0.0493 \\
        0.3498   &   5.5318   &   0.2882   &   0.4051   &   0.0219 \\ \hline \hline 
        Order   &   2.0429   &   2.0336   &   3.3843   &   3.299 \\ \hline 
    \end{tabular}
    \caption{Convergence rates for the two-dimensional Taylor-Green vortex at  $t_f = 0.5$ with $\nu=10^{-2}$.}
    \label{tab:tgv2d_10-2}
\end{table}

\begin{table}[htbp]
    \centering
    \begin{tabular}{|c || c |c || c |c |}
        \hline
        $h_{\text{max}}$ & $\P_1$ $\norm{e_u}_{L_1}$ & $\P_1$ $\norm{e_u}_{L_2}$ & 
        $\P_2$ $\norm{e_u}_{L_1}$ & $\P_2$ $\norm{e_u}_{L_2}$    \\
        \hline \hline
        0.4856   &   12.6204   &   0.6362   &   1.5963   &   0.0844 \\
        0.4466   &   8.8993   &   0.4525   &   1.1787   &   0.0639 \\
        0.3498   &   5.7806   &   0.3057   &   0.641   &   0.0357 \\  \hline  \hline 
        Order   &   2.2422   &   2.0922   &   2.7169   &   2.5713 \\  \hline 
    \end{tabular}\vspace*{0.2cm}
    \begin{tabular}{|c || c |c || c |c |}
        \hline
        $h_{\text{max}}$ & $\P_1$ $\norm{e_v}_{L_1}$ & $\P_1$ $\norm{e_v}_{L_2}$ & 
        $\P_2$ $\norm{e_v}_{L_1}$ & $\P_2$ $\norm{e_v}_{L_2}$  \\
        \hline \hline
        0.4856   &   12.1302   &   0.6146   &   1.5478   &   0.0814 \\
        0.4466   &   8.8654   &   0.4603   &   1.2019   &   0.067 \\
        0.3498   &   5.7598   &   0.3018   &   0.6182   &   0.0341 \\ \hline  \hline 
        Order   &   2.1569   &   2.0693   &   2.7805   &   2.6819 \\ \hline
    \end{tabular}
    \caption{Convergence rates for the two-dimensional Taylor--Green vortex at $t_f = 0.5$ with $\nu=10^{-5}$.}
    \label{tab:tgv2d_10-5}
\end{table}

\section{Conclusion} \label{sec:conclu}
In this work, we have presented the novel high-order ADER-AFE discontinuous Galerkin method for the solution of nonlinear hyperbolic systems on three-dimensional unstructured polyhedral meshes. The computational grids are created with an optimized polyhedral mesh generator that ensures high quality of each computational cell, making the mesh effective for numerical simulations. Using agglomerated continuous finite element basis functions on sub-tetrahedral grids within each polyhedral cell, the scheme achieves high spatial accuracy while maintaining computational efficiency through a quadrature-free formulation. The addition of an artificial viscosity shock capturing technique allows for robust handling of discontinuities which improves the scheme robustness. 
The numerical results demonstrate that the proposed method is accurate, robust and stable, by correctly capturing complex flow features and respecting the theoretical convergence orders. These promising results validate the potential of the ADER-AFE-DG approach to target more complex applications in fluid and solid mechanics.

Several perspectives can be proposed from this work. First, alternative limiting strategies could be explored to further improve the scheme's robustness and reduce numerical dissipation near shocks following the recent developments of monolithic sub-flux limiters \cite{vilar2025local} . Additionally, extending the scheme to an Arbitrary-Lagrangian-Eulerian (ALE) formulation will open new possibilities to simulate a broader class of problems, such as fluid-structure interaction and multi-material flows. Moreover, the quadrature-free nature of the method also makes it suitable for GPU acceleration and large-scale parallel simulations, potentially enabling applications in turbulent flows, aerospace engineering, and geophysical modeling. 


\bigskip

\bigskip


\acknowl{The Authors received financial support by Fondazione Cariplo and Fondazione CDP (Italy) under the project No. 2022-1895. WB also acknowledges funding by the \textit{Institute de Mathématiques pour la Planète Terre} (IMPT - France) under the project \textit{AAP2025-ERUPTA}.}

\funding{Fondazione Cariplo and Fondazione CDP (Italy) under the project No. 2022-1895.}

\dataava{Data will be made available on reasonable request.}

\complia{The authors certify that there is no actual or potential conflict of interest in relation to this article.}



\bigskip


\appendix

\section{Analytical definition of the basis functions and nodal coordinates for $\mathbb{P}_1$ to $\mathbb{P}_4$ on the reference tetrahedron} \label{app:basis_functions}

\begin{figure}[!htbp]
    \centering
    \includegraphics[trim=3.3cm 2cm 4cm 3cm, clip, width=0.48\textwidth]{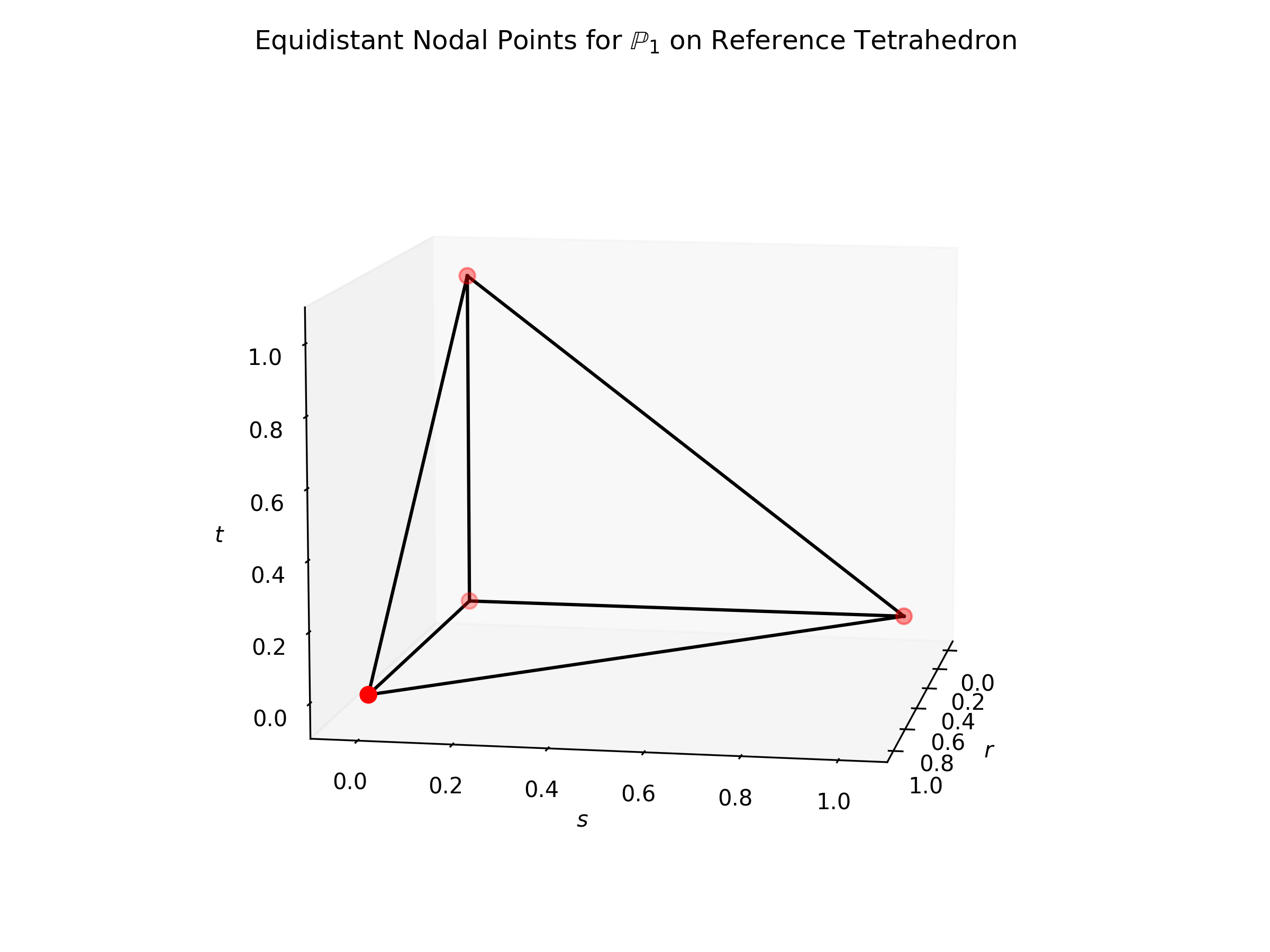}
    \includegraphics[trim=3.3cm 2cm 4cm 3cm, clip, width=0.48\textwidth]{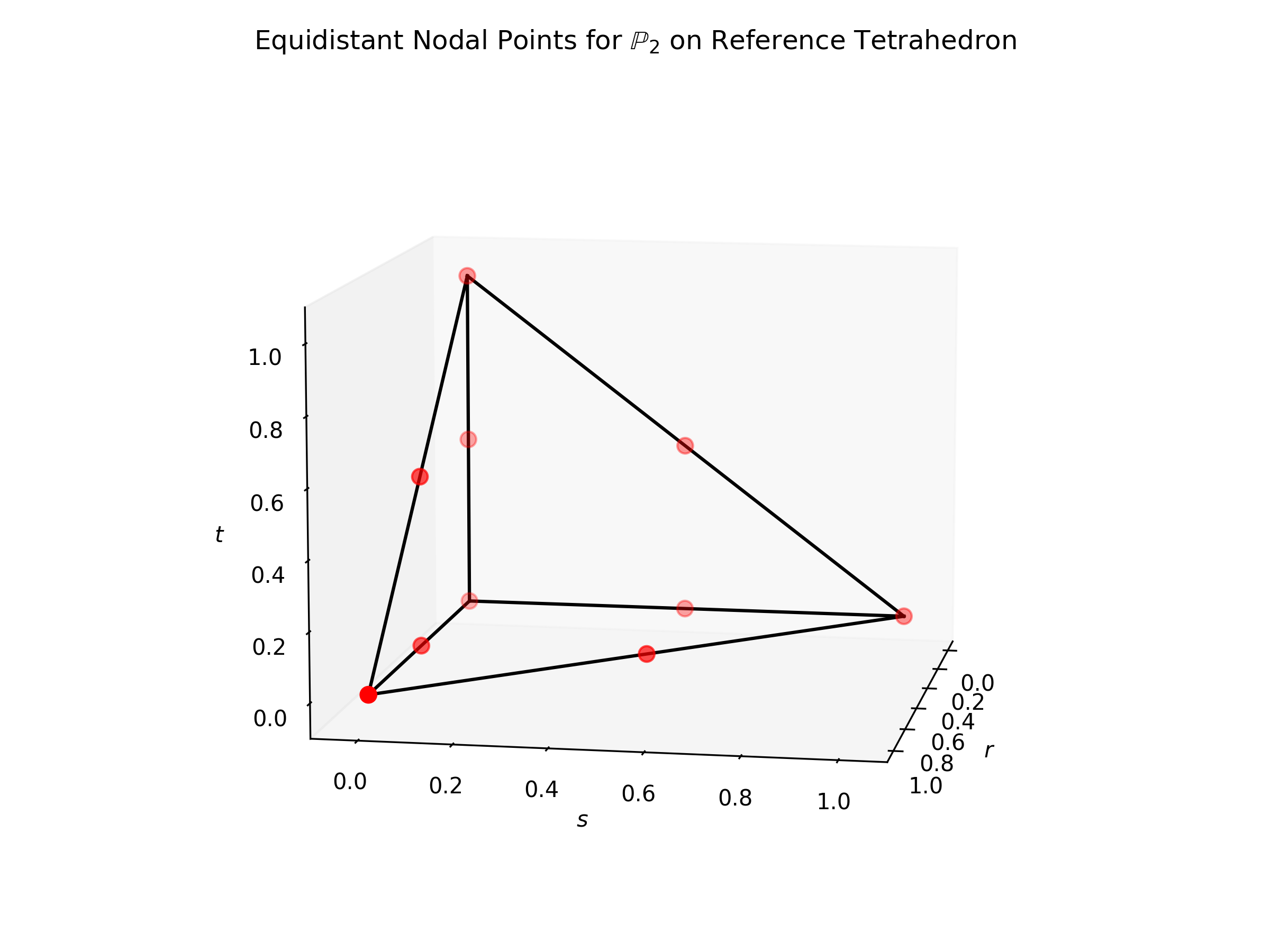}
    \includegraphics[trim=3.3cm 2cm 4cm 3cm, clip, width=0.48\textwidth]{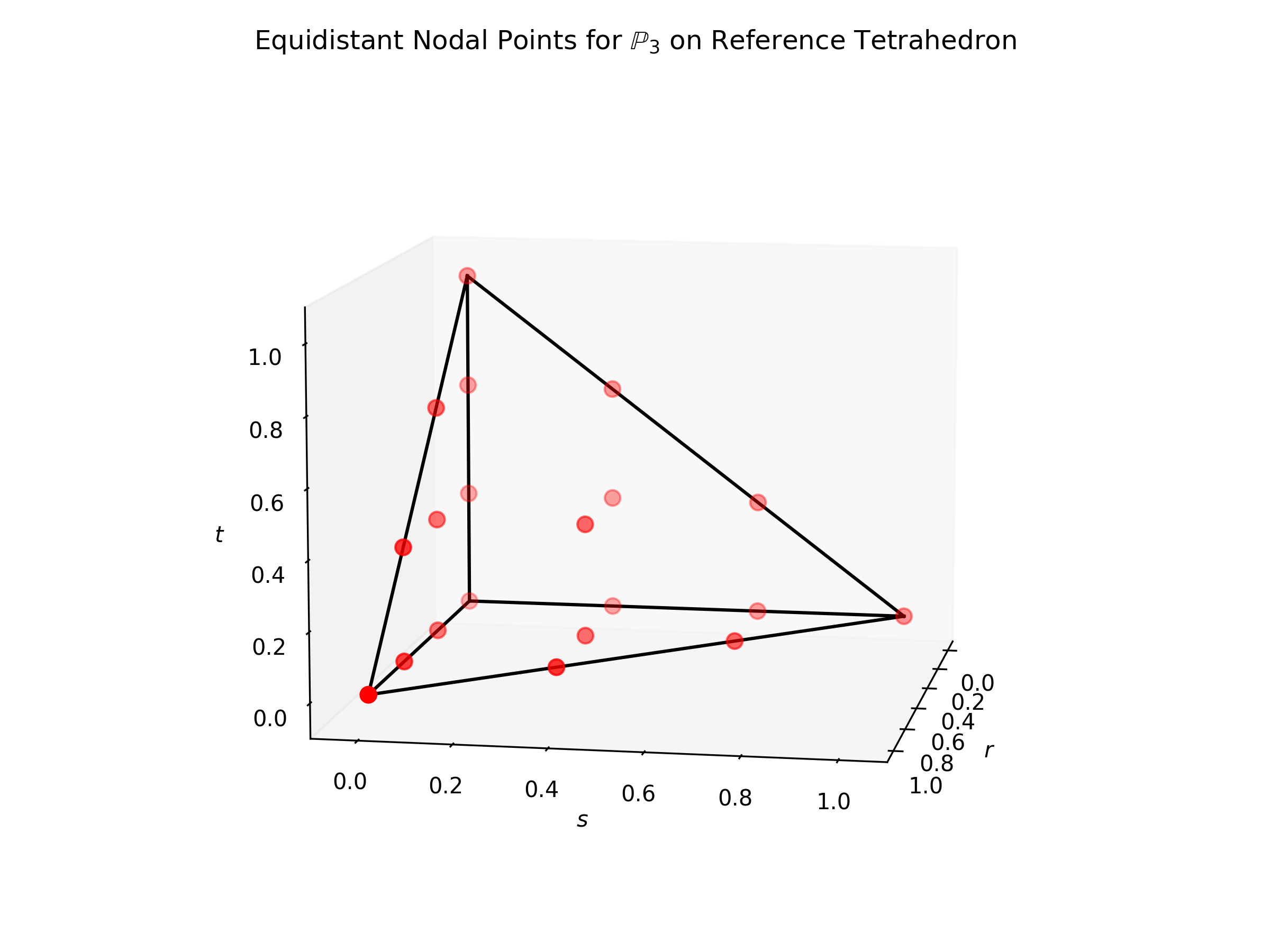}
    \includegraphics[trim=3.3cm 2cm 4cm 3cm, clip, width=0.48\textwidth]{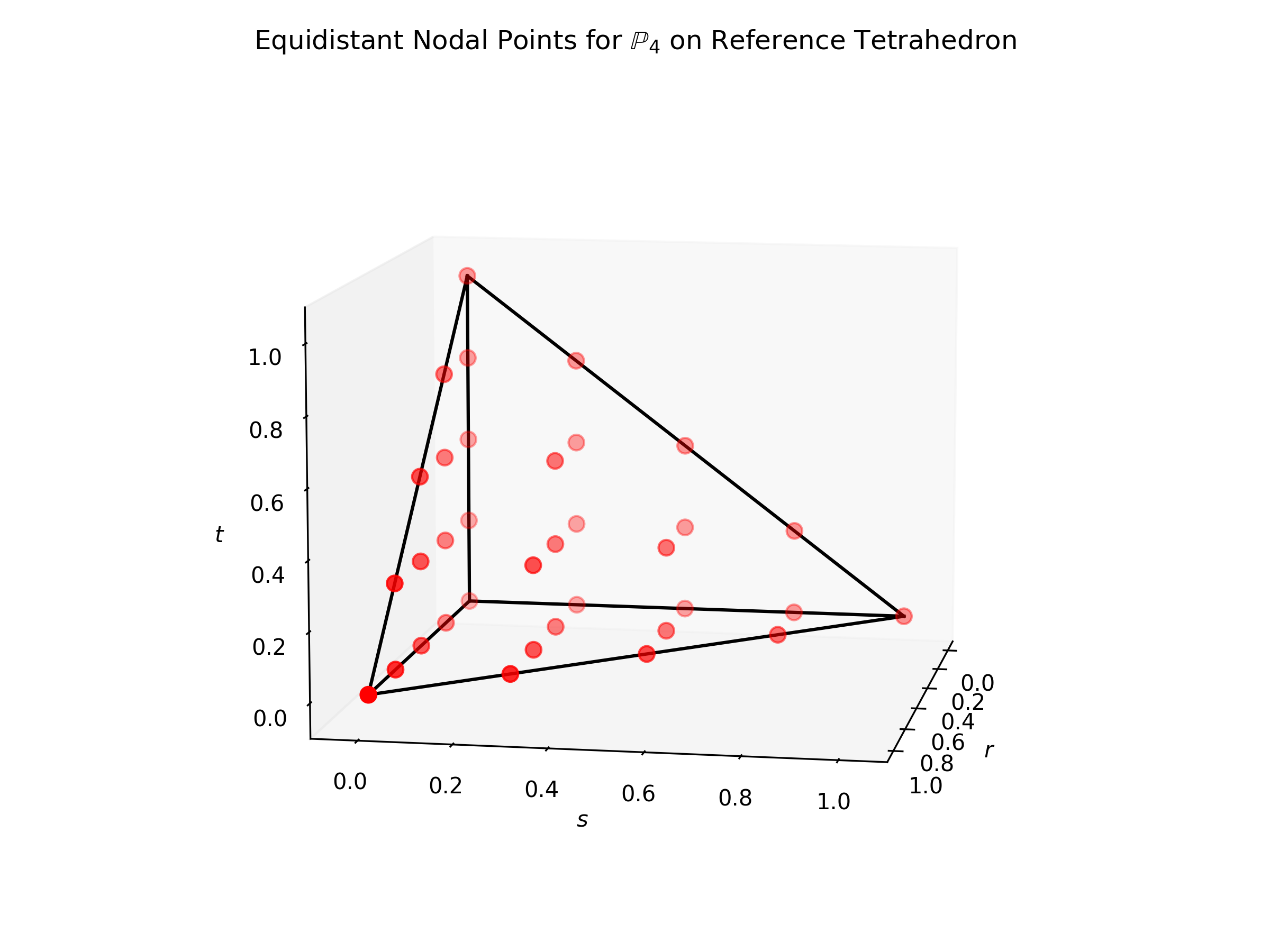}
    \caption{Elementary tetrahedral element: with a second (top left), third (top right), fourth (bottom left) and fifth (bottom right) order spatial discretization.}
    \label{fig:Tetra_all}
\end{figure}

In this appendix, we provide explicit expressions for the nodal Lagrange basis functions defined on the reference tetrahedron $\hat{T}$
\begin{equation*}
\hat{T} = \{ (r,s,t) \in \mathbb{R}^3 \mid r,s,t \geq 0, \quad r + s + t \leq 1 \},
\end{equation*}
with vertices at $(0,0,0), \, (1,0,0),\, (0,1,0)$ and $(0,0,1)$. Graphical illustrations of the corresponding elementary elements for several orders of accuracy are provided in Fig. \ref{fig:Tetra_all}.
The basis functions $\phi_i(r,s,t)$ satisfy the nodal property $\phi_i(\mathbf{x}_j) = \delta_{ij}$, where $\{\mathbf{x}_j\}$ are the interpolation nodes.

\subsection{Linear basis functions $\mathbb{P}_1$}

The linear basis functions correspond to the barycentric coordinates of the tetrahedron:
\begin{equation*}
\phi_1(r,s,t) = 1 - r - s - t, \quad \phi_2(r,s,t) = r, \quad \phi_3(r,s,t) = s, \quad \phi_4(r,s,t) = t.
\end{equation*}
The nodal points coordinates of Degree Of Freedom (DOF) are defined by:
\begin{equation*}
\begin{cases}
\mathbf{x}_1 = (0,0,0), \\
\mathbf{x}_2 = (1,0,0), \\
\mathbf{x}_3 = (0,1,0), \\
\mathbf{x}_4 = (0,0,1).
\end{cases}
\end{equation*}

\subsection{Quadratic basis functions $\mathbb{P}_2$}

The quadratic Lagrange basis functions are associated with the vertices and midpoints of the edges of $\hat{T}$:
\begin{equation*}
\begin{cases}
\phi_1 = \lambda_1 (2\lambda_1 - 1), \quad
\phi_2 = \lambda_2 (2\lambda_2 - 1), \quad
\phi_3 = \lambda_3 (2\lambda_3 - 1), \quad
\phi_4 = \lambda_4 (2\lambda_4 - 1), \\[1em]

\phi_5 = 4 \lambda_1 \lambda_2, \quad
\phi_6 = 4 \lambda_2 \lambda_3, \quad
\phi_7 = 4 \lambda_3 \lambda_1, \quad
\phi_8 = 4 \lambda_1 \lambda_4, \quad
\phi_9 = 4 \lambda_2 \lambda_4, \quad
\phi_{10} = 4 \lambda_3 \lambda_4.
\end{cases}
\end{equation*}
where $\lambda_1 = 1 - r - s - t$, $\lambda_2 = r$, $\lambda_3 = s$, $\lambda_4 = t$. The nodal points coordinates of DOF are defined by:
\begin{equation*}
\begin{cases}
\mathbf{x}_1 = (0,0,0), \\
\mathbf{x}_2 = (1,0,0), \\
\mathbf{x}_3 = (0,1,0), \\
\mathbf{x}_4 = (0,0,1), \\
\mathbf{x}_5 = \left(\frac{1}{2}, 0, 0\right), \quad
\mathbf{x}_6 = \left(\frac{1}{2}, \frac{1}{2}, 0\right), \quad
\mathbf{x}_7 = \left(0, \frac{1}{2}, 0\right), \\
\mathbf{x}_8 = \left(0, 0, \frac{1}{2}\right), \quad
\mathbf{x}_9 = \left(\frac{1}{2}, 0, \frac{1}{2}\right), \quad
\mathbf{x}_{10} = \left(0, \frac{1}{2}, \frac{1}{2}\right).
\end{cases}
\end{equation*}

\subsection{Cubic basis functions $\mathbb{P}_3$}

The cubic basis functions are constructed from barycentric coordinates and associated with nodes placed at vertices, edge points, face points, and interior points. 
The 20 DOF for $\mathbb{P}_3$ on $\hat{T}$ are:
\begin{equation*}
\begin{cases}
\text{Vertices:} \\
\mathbf{x}_1 = (0,0,0), \quad \mathbf{x}_2 = (1,0,0), \quad \mathbf{x}_3 = (0,1,0), \quad \mathbf{x}_4 = (0,0,1), \\
\text{Edge nodes (at } \frac{1}{3} \text{ and } \frac{2}{3} \text{ along each edge):} \\
\mathbf{x}_5 = \left(\frac{1}{3},0,0\right), \quad \mathbf{x}_6 = \left(\frac{2}{3},0,0\right), \\
\mathbf{x}_7 = \left(\frac{2}{3},\frac{1}{3},0\right), \quad \mathbf{x}_8 = \left(\frac{1}{3},\frac{2}{3},0\right), \\
\mathbf{x}_9 = \left(0,\frac{2}{3},0\right), \quad \mathbf{x}_{10} = \left(0,\frac{1}{3},0\right), \\
\mathbf{x}_{11} = \left(0,0,\frac{1}{3}\right), \quad \mathbf{x}_{12} = \left(0,0,\frac{2}{3}\right), \\
\mathbf{x}_{13} = \left(\frac{1}{3},0,\frac{2}{3}\right), \quad \mathbf{x}_{14} = \left(\frac{2}{3},0,\frac{1}{3}\right), \\
\mathbf{x}_{15} = \left(0,\frac{1}{3},\frac{2}{3}\right), \quad \mathbf{x}_{16} = \left(0,\frac{2}{3},\frac{1}{3}\right), \\
\text{Face interior node:} \\
\mathbf{x}_{17} = \left(\frac{1}{3},\frac{1}{3},0\right), \quad \mathbf{x}_{18} = \left(\frac{1}{3},0,\frac{1}{3}\right),   \\
\mathbf{x}_{19} = \left(0,\frac{1}{3},\frac{1}{3}\right), \quad \mathbf{x}_{20} = \left(\frac{1}{3},\frac{1}{3},\frac{1}{3}\right).
\end{cases}
\end{equation*}
And basis functions are defined by 
\begin{equation*}
\begin{cases}
\phi_1 = \frac{1}{2} \lambda_1 (3\lambda_1 - 1)(3\lambda_1 - 2), \quad
\phi_2 = \frac{1}{2} \lambda_2 (3\lambda_2 - 1)(3\lambda_2 - 2), \\
\phi_3 = \frac{1}{2} \lambda_3 (3\lambda_3 - 1)(3\lambda_3 - 2), \quad
\phi_4 = \frac{1}{2} \lambda_4 (3\lambda_4 - 1)(3\lambda_4 - 2), \\[1em]

\phi_5 = \frac{9}{2} \lambda_1 \lambda_2 (3\lambda_1 - 1), \quad
\phi_6 = \frac{9}{2} \lambda_1 \lambda_2 (3\lambda_2 - 1), \\

\phi_7 = \frac{9}{2} \lambda_1 \lambda_3 (3\lambda_1 - 1), \quad
\phi_8 = \frac{9}{2} \lambda_1 \lambda_3 (3\lambda_3 - 1), \\

\phi_9 = \frac{9}{2} \lambda_1 \lambda_4 (3\lambda_1 - 1), \quad
\phi_{10} = \frac{9}{2} \lambda_1 \lambda_4 (3\lambda_4 - 1), \\

\phi_{11} = \frac{9}{2} \lambda_2 \lambda_3 (3\lambda_2 - 1), \quad
\phi_{12} = \frac{9}{2} \lambda_2 \lambda_3 (3\lambda_3 - 1), \\

\phi_{13} = \frac{9}{2} \lambda_2 \lambda_4 (3\lambda_2 - 1), \quad
\phi_{14} = \frac{9}{2} \lambda_2 \lambda_4 (3\lambda_4 - 1), \\

\phi_{15} = \frac{9}{2} \lambda_3 \lambda_4 (3\lambda_3 - 1), \quad
\phi_{16} = \frac{9}{2} \lambda_3 \lambda_4 (3\lambda_4 - 1), \\[1em]

\phi_{17} = 27 \lambda_1 \lambda_2 \lambda_3, \quad
\phi_{18} = 27 \lambda_1 \lambda_2 \lambda_4, \\
\phi_{19} = 27 \lambda_1 \lambda_3 \lambda_4, \quad
\phi_{20} = 27 \lambda_2 \lambda_3 \lambda_4.
\end{cases}
\end{equation*}
where $\lambda_1 = 1 - r - s - t$, $\lambda_2 = r$, $\lambda_3 = s$, $\lambda_4 = t$.

\subsection{Quartic basis functions $\mathbb{P}_4$}

For the quartic basis functions, 35 nodal points are distributed equispaced inside the tetrahedron on vertices, edges, faces, and interior. The nodal points can be defined by all tuples $(r,s,t)$ with
\begin{equation*}
r,s,t \in \left\{0, \frac{1}{4}, \frac{1}{2}, \frac{3}{4}, 1 \right\}, \quad r + s + t \leq 1,
\end{equation*}
resulting in 35 points including:
\begin{itemize}
    \item[-] 4 vertices at $(0,0,0), (1,0,0), (0,1,0), (0,0,1)$,
    \item[-] 12 points on edges (at $\frac{1}{4}, \frac{1}{2}, \frac{3}{4}$),
    \item[-] 18 points on faces,
    \item[-] 1 interior point at $\left(\frac{1}{4}, \frac{1}{4}, \frac{1}{4}\right)$.
\end{itemize}
Explicit listing of all 35 points is lengthy but follows the same pattern. \\



\bibliographystyle{plain}
\bibliography{bibliography}

\end{document}